\documentclass[10pt]{amsart}
\usepackage{amscd,amssymb}
\usepackage{enumerate}
\usepackage[all]{xypic}
\usepackage{multirow}

\theoremstyle{plain}
\newtheorem{Thm}[equation]{Theorem}
\newtheorem{Cor}[equation]{Corollary}
\newtheorem{Prop}[equation]{Proposition}
\newtheorem{Lem}[equation]{Lemma}
\numberwithin{equation}{section}

\newcommand{\Irr}{\operatorname{Irr}}
\newcommand{\Hom}{\operatorname{Hom}}

\newcommand{\C}{\mathbb C}
\newcommand{\A}{\mathbb{A}}

\newcommand{\N}{\mathbb{N}}

\newcommand{\bm}{\begin{multline*}}
\newcommand{\tu}{\end  {multline*}}

\newcommand{\hence}{\Rightarrow}

\newcommand{\ie}{{\em i.e. }}

\title[On Shalika Periods]{On Shalika Periods and a Theorem of Jacquet-Martin}
\author{Wee Teck Gan and Shuichiro Takeda}
\address{Mathematics Department, University of California, San Diego, 9500 Gilman Drive, La Jolla,
92093} \email{wgan@math.ucsd.edu} \email{shtakeda@math.ucsd.edu}

\begin{document}

\begin{abstract}
Let $\pi$ be a cuspidal automorphic representation of $GL_4(\A)$
with central character $\mu^2$. It is known that $\pi$ has Shalika
period with respect to $\mu$ if and only if the $L$-function $L^S(s, \pi,
\bigwedge^2\otimes\mu^{-1})$ has a pole at $s=1$.
In \cite{JM}, Jacquet and Martin considered the analogous question for cuspidal
representations $\pi_D$ of the inner form $GL_2(D)(\A)$, and
obtained a partial result via the relative trace formula. In this
paper, we provide a complete solution to this problem via the method of
theta correspondence, and give necessary and sufficient
conditions for the existence of Shalika period for $\pi_D$. We also
resolve the analogous question in the local setting.
\end{abstract}

\maketitle

 \section{\bf Introduction}
Let $F$ be a number field with adele ring $\A$, and let $D$ be a (possibly split)
quaternion algebra over $F$. We consider the linear
algebraic group $GL_2(D)$, so that if $D$ is split, then $GL_2(D)
\cong GL_4$. The group $GL_2(D)$ is thus an inner form of $GL_4$. Let
$\pi_D$ be a cuspidal automorphic representation of $GL_2(D)(\A)$ and
assume that its central character $\omega_{\pi_D}$ is a square, say
\[  \omega_{\pi_D} = \mu^2. \]
One may consider the Shalika period  of $\pi_D$ with respect to
$\mu$. More precisely, $GL_2(D)$ has a parabolic $F$-subgroup
$P_D=M_D\cdot N_D$ with Levi factor and unipotent radical given by:
 \[  \begin{cases}
 M_D \cong D^{\times} \times D^{\times}; \\
N_D \cong D. \end{cases} \] Let $\psi_D$ be the nondegenerate
unitary character of $N_D(F) \backslash N_D(\A)$ defined by
$\psi_D(x) = \psi(Tr_D(x))$ for $x \in N_D(\A) \cong D_{\A}$. Its
stabilizer in $P_D$ is the Shalika subgroup
\[  \tilde{S}_D = \Delta D^{\times} \cdot N_D \]
and we may extend $\psi_D$ to a character of $\tilde{S}_D$ via:
\[  \psi_D(h\cdot n) = \mu(\mathbb{N}_D(h)) \cdot \psi_D(n) \quad \text{for $h \in D^{\times}$ and $n \in N_D$.} \]
We shall in fact mostly be concerned with the quotient group
\[  S_D = \tilde{S}_D/ \Delta \mathbb{G}_m = PD^{\times} \cdot N_D. \]
The Shalika period of $\pi_D$ is the linear form on $\pi_D$ defined
by: \vskip 5pt
\[ \mathcal{S}_D:  f \mapsto  \int_{S_D(F) \backslash S_D(\A)}
f(nh) \cdot \mu(\mathbb{N}_D(h))^{-1} \cdot  \psi_D(n)^{-1} \, dn \,
dh. \] \vskip 5pt \noindent We say that $\pi_D$ has Shalika period
with respect to $\mu$ if the linear form $\mathcal{S}_D$ is
non-zero. If $\mu$ is trivial, then we simply say that $\pi_D$ has
Shalika period. In the following, if $D$ is split, we shall suppress
the symbol $D$ from the above notations. So, for example, $GL_4$
has a parabolic subgroup $P = M \cdot N$. \vskip 10pt

Now suppose that $D$ is split, so that $GL_2(D)=GL_4$. There is a well-known
theorem of Jacquet and Shalika \cite{JS} that relates the
existence of Shalika period on $GL_{2n}$ to the existence of poles of a twisted exterior square
$L$-function. For the case of $GL_4$, their theorem reads as follows.

\begin{Thm}[Jacquet-Shalika]\label{T:JS}
Let $\pi$ be a cuspidal automorphic representation of $GL_4(\A)$
whose central character is $\mu^2$. Then the following are
equivalent:
\begin{enumerate}[(a)]
\item $\pi$ has Shalika period with respect to $\mu$.
\item The (incomplete) twisted exterior square $L$-function $L^S(s, \pi, \bigwedge^2\otimes\mu^{-1})$
        has a pole at $s= 1$.
\end{enumerate}
\end{Thm}

Now it is natural to ask whether the same
theorem holds when $D$ is not split. In
their recent paper \cite{JM}, Jacquet and Martin obtained the following result
by using the relative trace formula.

\begin{Thm}[Jacquet-Martin]  \label{T:JM}
Suppose that $D$ is a quaternion division algebra and $\pi_D$ is a
cuspidal representation of $GL_2(D)(\A)$ which has a cuspidal
Jacquet-Langlands lift $\pi$ to $GL_4$. Further assume that \vskip
5pt \noindent (i) $D$ is non-split at some archimedean place; \vskip
5pt \noindent (ii) $\pi_D$ has trivial central character; \vskip 5pt
\noindent (ii) $\pi_D$  has at least one supercuspidal local
component at a place where $D$ splits. \vskip 5pt \noindent  Then
\[  \text{$\pi_D$ has Shalika period $\Longrightarrow$ $\pi$ has Shalika period.} \]
\end{Thm}
\vskip 5pt

In light of Theorem 1.1, their theorem shows that,
modulo some technicalities, if $\pi_D$ has Shalika period, then the
exterior square $L$-function $L^S(s, \pi_D,\bigwedge^2)$ has a pole
at $s=1$.  Then the natural question to ask is  whether the converse is also true; in other words, whether the analog of Thm. \ref{T:JS} remains true for non-split $D$.\\

In this paper, we  resolve this question completely by giving a characterization of the existence of Shalika period for all cuspidal representations $\pi_D$. It turns out that the converse of the Jacquet-Martin theorem does not hold and must be augmented with a certain local condition. Our main theorem is:

\begin{Thm}\label{T:main}
Let $D$ be a quaternion division algebra and $\Sigma_D$ the set of
places at which $D$ ramifies. Further assume that $D$ splits at
every archimedean place.
\vskip 10pt

\noindent (i) Suppose that $\pi_D$ is a cuspidal
automorphic representation on $GL_2(D)(\A)$ with central character
$\mu^2$ and whose Jacquet-Langlands lift $JL_(\pi_D)$ to
$GL_4(\A)$ is cuspidal. Then the following are equivalent.
\vskip 5pt

\begin{enumerate}[(A)]
\item $\pi_D$ has Shalika period with respect to $\mu$.
\vskip 5pt

\item The (incomplete) twisted
exterior square $L$-function $L^S(s, \pi_D, \bigwedge^2\otimes\mu^{-1})$
       has a pole at $s= 1$, and for all $v \in \Sigma_D$, $\pi_{D,v}$ is not of the form
        $Ind^{GL_2(D_v)}_{P_{D,v}} \delta_{P_D}^{1/2} \cdot (\tau_{D,1,v} \boxtimes \tau_{D,2,v})$
        where  $\tau_{D,i,v}$ are representations of $D_v^{\times}$ with central character $\mu_v$.
\end{enumerate}
\vskip 10pt

\noindent (ii) Suppose that $\pi_D$ is cuspidal with central character $\mu^2$ but
its Jacquet-Langlands lift $JL(\pi_D)$ to $GL_4(\A)$ is not cuspidal. In this case,  $JL(\pi_D)$ is contained in the residual spectrum and is isomorphic to the unique irreducible quotient of
$Ind^{GL_4}_P \delta_P^{1/2} \cdot \left( \tau |-|^{1/2} \boxtimes \tau |-|^{-1/2} \right)$ for a cuspidal
representation $\tau$ of $GL_2(\A)$. Then the following are equivalent:
\vskip 5pt

\noindent (C)  $\pi_D$ has Shalika period with respect to $\mu$.
\vskip 5pt

\noindent (D) $\mu$ is equal to the central character $\omega_{\tau}$ of $\tau$.
\vskip 5pt

\noindent (E)  The (incomplete) twisted
exterior square $L$-function $L^S(s, \pi_D, \bigwedge^2\otimes\mu^{-1})$ has a pole at $s= 2$.
 \end{Thm}
\vskip 10pt

The reason for assuming that $D$ is split at every archimedean place
in the above theorem is that we make use of recent results of
Badulescu [B] concerning the Jacquet-Langlands correspondence and
this assumption is present in his work. In fact, for the Jacquet-Martin theorem, which is part of the implication $(A) \Longrightarrow (B)$, one does not need the assumption that $D$ be split at all archimedean places.

\vskip 10pt

Clearly the interesting point in our theorem is the local condition
in (B), which is not present in the split case. Let us briefly
explain the origin of this local condition. For each place $v$ of
$F$ and a representation $\pi_{D,v}$ of $GL_2(D_v)$ with central
character $\mu_v^2$, we say that $\pi_{D,v}$ has local Shalika
period with respect to $\mu_v$ if
\[ \Hom_{\Delta D_v^{\times} \cdot N_{D, v}} (\pi_{D,v},  (\mu \circ \N_{D,v}) \boxtimes \psi_{D_v}) \ne 0. \]
\vskip 5pt \noindent  It is known that this Hom space has dimension
at most 1. One may consider the problem of existence of local
Shalika periods, and indeed we will show that a generic
representation $\pi_v$ of $GL_4(F_v)$ has a local Shalika period if
and only if its Langlands parameter factors through $GSp_4(\C)$, \ie
is of symplectic type. Suppose that this holds and $\pi_{D,v}$ is
the local Jacquet-Langlands lift of $\pi_v$ to $PGL_2(D)$. Then it
is possible that $\pi_{D,v}$ does not have local Shalika period.
Indeed, whether $\pi_{D,v}$ has a local Shalika period or not is an
issue addressed by a special case of the local Gross-Prasad
conjecture, and thus it is controlled by a local epsilon factor
condition. To show the implication (B)$\hence$(A), one needs (at
least) these local epsilon factor conditions to be satisfied.
\vskip 10pt

However, let us mention here that even if the local epsilon factor
conditions are satisfied, it turns out that they are not sufficient
for the global representation $\pi_D$ to have Shalika period. In
fact, we prove the following perhaps somewhat surprising result:
\vskip 5pt

\begin{Thm}\label{T:counter}
Assume that $D$ is split at every archimedean place. There are
cuspidal representations $\pi_D$ of $PGL_2(D)$, with a cuspidal
Jacquet-Langlands lift to $PGL_4$, satisfying: \vskip 5pt
\begin{itemize}
\item[(i)] for all places $v$, $\pi_{D,v}$ has local Shalika period, and
\vskip 5pt

\item[(ii)] $L^S(s, \pi_D, \bigwedge^2)$ has a pole at $s= 1$, but
\vskip 5pt

\item[(iii)]  $\pi_D$ does not have global Shalika period.
\end{itemize}
\end{Thm}
\vskip 10pt

We should mention that though the study of the existence of local
Shalika periods elucidates the nature of the local conditions in our
theorem, the proofs of the above global theorems are largely independent of
this local study. The exception is Prop. \ref{P:bad-type}, whose proof relies on
the local study of Section \ref{S:explicit}.
\vskip 10pt

Our main local results, which complete some initial work of D. Prasad, are summarized as follows:
\vskip 5pt

\begin{Thm}
Let $F_v$ be a non-archimedean local field and $D_v$ the unique
quaternion division algebra over $F_v$. \vskip 5pt

\noindent  (i) Let $\pi_v$ be a generic representation of
$GL_4(F_v)$ with central character $\mu_v^2$. Then $\pi_v$ has
Shalika period with respect to $\mu_v$ if and only if  the Langlands
parameter $\varphi_{\pi_v}$ of $\pi_v$ factors through $GSp_4(\C)$ with similitude
character $\mu_v$.
\vskip 5pt

 If $\pi_v$ is a discrete series representation,
the above conditions are equivalent to:
\[  \text{$L(s,\pi_v, \bigwedge^2 \otimes \mu_v^{-1})$  has a pole at $s=0$.} \]
 Here,  the local L-function $L(s,\pi_v, \bigwedge^2
\otimes \mu_v^{-1})$ is that defined by Shahidi, and is equal to the Artin L-function
$L(s, \bigwedge^2 \varphi_{\pi_v} \otimes \mu_v^{-1})$ by a recent result of Henniart [He].
  \vskip 10pt

\noindent (ii) Let $\pi_{D,v}$ be a representation of $GL_2(D_v)$
with generic Jacquet-Langlands transfer $\pi_v=JL(\pi_{D,v})$ on
$GL_4$ and central character $\mu_v^2$. Then $\pi_{D,v}$ has Shalika
period with respect to $\mu_v$ if and only if
\[ \epsilon(1/2, (\bigwedge^2 \varphi_{\pi_v} \otimes \mu_v^{-1}) \otimes S_2)= -1. \]
Here $S_2$ denotes the 2-dimensional representation of $SL_2(\C)$ (which is the Langlands parameter of the Steinberg representation of $GL_2$) and
$(\bigwedge^2 \varphi_{\pi_v} \otimes \mu_v^{-1}) \otimes S_2$ is a representation of the Weil-Deligne group $W_{F_v} \times SL_2(\C)$ of $F_v$.
Moreover, if the above holds, then
the Jacquet-Langlands transfer $\pi_v$ has Shalika period with
respect to $\mu_v$.
\end{Thm}
\vskip 5pt

In fact, we also determine whether a generalized Speh representation of
$GL_4$ or $GL_2(D)$ has Shalika period with respect to $\mu$; the result is contained in
Thm. \ref{T:speh}.  Since those local results are largely independent of our global
results, we will take them up at the end of the paper (Section
\ref{S:local} and \ref{S:explicit}).\\

\vskip 10pt

The main technique used in this paper is the theta
correspondence (for similitudes). Indeed, consider the quadratic
space
\[  (V_D, q_D)  = (D , \mathbb{N}_D) \oplus \mathbb{H} \]
where $\mathbb{H}$ is the hyperbolic plane. Then one has
\[  GSO(V_D) \cong (GL_2(D)  \times GL_1)/ \{ (z, z^{-2}): z \in GL_1\}. \]
To see this, note that the quadratic space $V_D$ can also be
described as the space of $2 \times 2$-Hermitian matrices with
entries in $D$, so that a typical element has the form
\[  (a,d; x) = \left( \begin{array}{cc}
a & x \\
\overline{x} & d \end{array} \right),  \qquad \text{$a, d \in F$ and
$x \in D$},  \] equipped with the quadratic form
\[\  -\det:  \left( \begin{array}{cc}
a & x \\
\overline{x} & d \end{array} \right) \mapsto -ad + \mathbb{N}_D(x).
\] The action of $GL_2(D)  \times GL_1$ on this space is given by
\vskip 5pt
\[  (g,z)(X) = z \cdot g \cdot X \cdot \overline{g}^t \]
\vskip 5pt \noindent Observe that an irreducible representation of
$GSO(V_D)$ is of the form $\pi \boxtimes \mu$ where $\pi$ is a
representation of $GL_2(D)$ and $\mu$ is a square root of the
central character of $\pi$. This is precisely the data needed to
define a Shalika period.
 \vskip 10pt

One can thus consider the theta correspondence for the (almost) dual
pair
\[
    GSp_4 \times GSO(V_D).
\]
When $D$ is split, this theta correspondence can be used to prove
the (weak) lifting of globally generic cuspidal representations of
$GSp_4$ to $GL_4$; this is a well-known result of Jacquet,
Piatetski-Shapiro and Shalika that was announced almost thirty
years ago, but whose proof was never published. However, most of the
details of their proof can be found in the paper \cite{So}, where
Soudry made use of this same dual pair to prove the strong
multiplicity one theorem for globally generic cuspidal
representations of $GSp_4$. Locally, a preliminary study of this
theta correspondence has been conducted by Waldspurger as {\em un
exercice} \cite{W}.
\vskip 10pt

The understanding of this theta correspondence for both split and
non-split $D$ underlies the results of this paper. In addition, for
the proof of Thm. \ref{T:main} and Thm. \ref{T:counter}, a key tool is a Rankin-Selberg
integral representation of the degree 5 L-function of a cuspidal
representation of $Sp_4$, which was first discovered by Andrianov in
the classical setting and recast in the adelic setting by
Piatetski-Shapiro and Rallis [PSR].

\vskip 10pt

Finally, we should mention that the theorem of Jacquet and Martin
has an obvious analog for $D$ an arbitrary division algebra of
degree $d$. Though the use of theta correspondence gives a simple
proof in the case when $D$ is quaternion, it has no hope of
addressing the general case. On the other hand, one fully expects
the relative trace formula approach of Jacquet and Martin to work
for general $d$, as long as one can master the analytic
difficulties.  However, as we learned shortly after the completion of this paper,
the analog of the Jacquet-Martin theorem (i.e. Thm. \ref{T:JM}) for general $D$
has been proven in a recent paper of Jiang-Nien-Qin [JNQ] using entirely different methods.
Our Theorem \ref{T:main} indicates that the correct converse statement in the case of general $D$ may be quite delicate.

 \vskip 15pt

\noindent{\bf Acknowledgments:} We have benefitted from many
illuminating email correspondences with Dipendra Prasad, as well as
from his various papers on the local Shalika model. We take this
opportunity to thank him for his help and for his many comments on
an early draft of this paper. We also thank Kimball Martin for
discussions concerning his paper with Jacquet, Ioan Badulescu for
conversations related to his work on the Jacquet-Langlands
correspondence, and last but not least, Gordan Savin for a catalytic
conversation which led us to work on this problem.  W.T. Gan's
research is partially supported by NSF grant DMS-0500781.

\vskip 15pt

\section{\bf Theta Correspondence for Similitudes}  \label{S:theta}

In this section, we give a brief introduction of the necessary
background on theta correspondence for similitudes. We shall follow
the reference \cite{Ro2} closely.
 \vskip 5pt

Let us begin by establishing some more group theoretic notations.
First, let us fix the isomorphism
\[  GSO(D, \mathbb{N}_D) \cong (D^{\times} \times D^{\times})/ \{(z,z^{-1}): z \in GL_1\} \]
via the action of the latter on $D$ given by
\[  (\alpha, \beta) \mapsto  \alpha x \overline{\beta}. \]
In particular, if $D$ is split, then $GSO(D, \mathbb{N}_D)$ is the
split orthogonal group $GSO(2,2)$ and we have:
\[  GSO(2,2) \cong (GL_2 \times GL_2)/ \{(z,z^{-1}): z \in GL_1\} \]
In any case, an irreducible representation of $GSO(D, \mathbb{N}_D)$ is thus of the form $\tau_{D,1} \boxtimes \tau_{D,2}$ where the central characters of $\tau_{D,1}$ and $\tau_{D,2}$ are equal.
\vskip 5pt

Also, as we explained in the introduction, there is a natural
isomorphism
\[  (GL_2(D) \times GL_1)/ \{(z,z^{-2}): z \in GL_1\} \overset{\sim}{\rightarrow} GSO(V_D). \]
In particular if $D$ is split, we simply write $V$ for $V_D$ and
\[  (GL_4 \times GL_1)/ \{(z,z^{-2}): z \in GL_1\} \overset{\sim}{\rightarrow} GSO(V). \]
The similitude factor of $GSO(V_D)$ is: \vskip 5pt
\[  \lambda_{D}: (g,z) \mapsto N(g) \cdot z^2, \]
\vskip 5pt \noindent where $N$ is the reduced norm on the central
simple algebra $M_2(D)$. Thus,
 \[  SO(V_D) = \{ (g,z) \in GSO(V_D): N(g) \cdot z^2 = 1\}. \]
\vskip 5pt

Under the above isomorphism, the parabolic subgroup $P_D \times
GL_1$ is identified with the stabilizer of the line in $V_D$ spanned
by $(1,0; 0_D)$.  Indeed, using matrix representation with respect
to the decomposition $V_D = F \cdot (1,0) \oplus D \oplus F \cdot
(0,1)$, we have: \vskip 5pt
\[  \left( \left( \begin{array}{cc}
\alpha &  0 \\
0 & \beta \end{array} \right), z \right)  \mapsto \left(
\begin{array}{ccc}
z \cdot \N(\alpha) & & \\
& z \cdot (\alpha, \beta) & \\
&  &  z \cdot \N(\beta) \end{array} \right) \] and
\[
\left( \begin{array}{cc}
1 &  y \\
0 & 1 \end{array} \right)  \mapsto \left( \begin{array}{ccc}
1 & Tr_D(\overline{y} - ),     & \N(y)  \\
& 1 & y \\
&  &  1  \end{array} \right). \] \vskip 5pt

Observe that there is an embedding
\[  \iota: S_D \cong PD^{\times} \cdot N_D \hookrightarrow GSO(V_D). \]
The embedding of $N_D$ is the one given above, whereas $PD^{\times}$
is embedded via:
\[ \iota(\alpha)=
\left( \left( \begin{array}{cc}
\alpha &  \\
  & \alpha \end{array} \right) ,  \N \alpha^{-1} \right) \in (GL_2(D) \times GL_1)/\{ (z,z^{-2}) \}. \]
Thus, if $\pi \boxtimes \mu$ is a cuspidal representation of
$GSO(V_D)$, then the Shalika period on $\pi$ with respect to $\mu$
is simply the linear form on $\pi \boxtimes \mu$ given by
\[  f \mapsto  \int_{S_D(F) \backslash S_D(\A)}
f( \iota(n)) \cdot \overline{\psi_D(\iota(n))} \, dn \] where
$\psi_D$ is extended from $N_D$ to $S_D$ by requiring that $\psi_D$
be trivial on $PD^{\times}$.

\vskip 10pt

Now let $W$ be the 4-dimensional symplectic vector space
and fix a Witt decomposition $W = X \oplus Y$.  Let $P(Y) = GL(Y)
\cdot N(Y)$ be the parabolic subgroup stabilizing the maximal
isotropic subspace $Y$. Then
\[  N(Y) = \{ b \in Hom(X,Y) : b^t  = b \}, \]
where $b^t \in Hom(Y^*, X^*) \cong Hom(X,Y)$. \vskip 5pt

Fix a unitary character $\psi$ of $F \backslash \A$ and consider the
Weil representation $\omega_D$ associated to $\psi$ for the dual
pair  $Sp(W)(\A) \times O(V_D)(\A)$. It can be realized on $S((X
\otimes V_D)(\A))$ and the action of $P(Y) \times O(V_D)$ is given
by the usual formulas: \vskip 5pt
\[  \begin{cases}
\omega_D(h)\phi(x) = \phi(h^{-1} x), \quad \text{for $h \in O(V_D)$;} \\
\omega_D(a)\phi(x) = |\det_Y(a)|^{\frac{1}{2} \dim V_D} \cdot \phi(a^{-1} \cdot x), \quad \text{for $a \in GL(Y)$;} \\
\omega_D(b) \phi(x) = \psi( \langle bx, x \rangle) \cdot \phi(x),
\quad \text{for $b \in N(Y)$,} \end{cases}
\]
\vskip 5pt \noindent where $\langle -, -\rangle$ is the natural
symplectic form on $W \otimes V_D$. To describe the full action of
$Sp(W)$, one needs to specify the action of a Weyl group element,
which acts by a Fourier transform.

\vskip 10pt

Now let
\[  R_D =\{ (g,h)  \in GSp(W) \times GO(V_D): \lambda(g) \cdot \lambda_D(h)=1 \},  \]
where the $\lambda$'s refer to the similitude factor of the relevant
group. Note that this differs from the normalization in \cite{Ro2}.
The Weil representation can then be extended in a natural way to the
group $R_D(\A)$, via:
\[  \omega_D(g,h)\phi = |\lambda_D(h)|^{- \frac{1}{2}\dim V_D} \omega_D(g_1, 1)(\phi \circ h^{-1}) \]
where
\[  g_1 = g \left(  \begin{array}{cc}
\lambda(g)^{-1} & 0 \\
0 & 1  \end{array} \right) \in Sp(W). \] Observe that the central
elements $(t,t^{-1}) \in R_D$ act trivially. We shall in fact only
be interested in the action of
\[  R_D^0 = \{ (g,h) \in R_D: h \in GSO(V_D) \}. \]

\vskip 10pt

For $\phi \in \mathcal{S}((X \otimes V_D)(\A))$ and $(g,h) \in
R_D(\A)$, set
\[  \theta(\phi)(g,h)  = \sum_{x \in (X \otimes V_D)(F)} \omega_D(g,h)\phi(x). \]
Then $\theta(\phi)$ is a function of moderate growth on $R_D(F)
\backslash R_D(\A)$. If $\pi_D \boxtimes \mu$ is a cuspidal
representation of $GSO(V_D)$ and $f \in \pi_D\boxtimes \mu$, we set
\[  \theta(\phi,f)(g) = \int_{SO(V_D)(F)\backslash SO(V_D)(\A)} \theta(\phi)(g,h_1h) \cdot \overline{f(h)} \, dh \]
where $h_1$ is any element of $GSO(V_D)$ such that $\lambda_D(h_1) =
\lambda(g)$. Moreover, set
\[  \Theta(\pi_D\boxtimes \mu) = \langle \theta(\phi,f): \phi \in \omega_D, f \in \pi_D \boxtimes \mu \rangle. \]
Then $\theta(\phi,f)$ is an automorphic form (possibly zero) on
$GSp(W)$ and $\Theta(\pi_D \boxtimes \mu)$ is an automorphic
representation (possibly zero) of $GSp(W)(\A)$ whose central character is
equal to that of $\pi_D$.  Similarly, starting from a cuspidal
representation $\sigma$ of $GSp(W)$, we have the automorphic
representation $\Theta_D(\sigma)$ of $GSO(V_D)$.
 \vskip 10pt

Many questions about these similitude theta liftings can be easily reduced to the analogous questions for the isometry case. We highlight two such questions here.
\vskip 10pt

The first such question is
the vanishing or non-vanishing of $\Theta(\pi_D \boxtimes \mu)$.
If $res$ denotes the restriction of functions from
a similitude group to the corresponding isometry group, then $res(\pi_D \boxtimes \mu)$ is a nonzero cuspidal representation of $SO(V_D)$ which is possibly reducible. From the definition of $\theta(\phi,f)$, it is immediate that
$\Theta(\pi_D \boxtimes \mu)$ is nonzero iff the global theta lift of $res(\pi_D \boxtimes \mu)$ from $SO(V_D)$ to $Sp(W)$ is nonzero. Similarly, if $\sigma$ is a cuspidal representation of $GSp(W)$, then $\Theta_D(\sigma)$ is nonzero iff the global theta lift of $res(\sigma)$ from $Sp(W)$ to $SO(V_D)$
is nonzero.
\vskip 10pt

The second such question is the cuspidality of $\Theta(\pi_D \boxtimes\mu)$. Again, it is evident from the definition that $\Theta(\pi_D \boxtimes \mu)$ is contained in the space of cusp forms of $GSp(W)$ iff
the global theta lift of $res(\pi_D \boxtimes \mu)$ is contained in the space of cusp forms of $Sp(W)$.
Now in the isometry case, if one has a tower of theta liftings in the sense of Rallis, then a standard result in the theory is the so-called tower property of theta correspondence.  This says
that the global theta lift of a cuspidal representation of an isometry group to a particular step in the tower
is contained in the space of cusp forms iff its theta lift to the previous step of the tower vanishes. Together with the above discussion, one sees immediately that the same statement applies to the theta liftings for similitude groups. Moreover, after the first nonzero lift, the theta lifts to higher steps of the tower do not vanish and are not contained in the space of cusp forms (though its intersection with the space of cusp forms may be nonzero because we are working with $SO$ rather than $O$).

\vskip 15pt
The following lemma will be used in this paper.
\begin{Lem}\label{L:generic_non-vanish}
Let $\sigma$ be a globally generic cuspidal representation of $GSp_4$. Then the global theta lift $\Theta(\sigma)$ of $\sigma$ to $GSO(V)=GSO(3,3)$ is globally generic and thus is nonzero.
\end{Lem}
\begin{proof}
This is essentially the main theorem of \cite{GRS}. There they considered the isometry groups, but by our discussion above, it is easy to see that their theorem applies to the similitude case.
\end{proof}

We now note: \vskip 5pt

\begin{Prop}  \label{P:JL-theta}
(i) Suppose that $\pi_D \boxtimes \mu$ is a cuspidal representation of
$GSO(V_D)$ such that the Jacquet-Langlands lift $JL(\pi_D)$ of
$\pi_D$ to $GL_4$ is cuspidal. If $\Theta(\pi_D\boxtimes \mu)$ is
non-zero, then $\Theta(\pi_D\boxtimes \mu)$ is contained in the
space of cusp forms of $GSp_4$.
\vskip 10pt

(ii) Suppose now that $JL(\pi_D)$ is non-cuspidal.  If $\Theta(\pi_D\boxtimes \mu)$ is
non-zero, then $\Theta(\pi_D\boxtimes \mu)$ does not contain any globally generic cuspidal representation  of $GSp_4$.
\end{Prop}
\vskip 5pt

\begin{proof}
(i) By the tower propery of theta correspondence, if
$\Theta(\pi_D\boxtimes \mu)$ is non-cuspidal, then the theta lift of
$\pi_D\boxtimes \mu$ to $GL_2$ (which is the lower step of the
tower) is nonzero cuspidal.
Denote this cuspidal representation of
$GL_2$ by $\Sigma$. Consider the theta lift of $\Sigma$ to $GSO(V) =
GSO(3,3)$. Since the theta lift of $\Sigma$ to $GSO(2,2)$ is
well-known to be nonzero, it follows by the tower property again that its theta lift
$\Theta(\Sigma)$ to $GSO(V)$ is nonzero and not contained in the
space of cusp forms. But any irreducible subquotient of
$\Theta(\Sigma)$ is nearly equivalent to the cuspidal
$JL(\pi_D)\boxtimes \mu$. This contradicts the generalized strong
multiplicity one theorem of Jacquet-Shalika.
\vskip 10pt

(ii) If  $\Theta(\pi_D\boxtimes \mu)$ contains a globally generic
cuspidal representation $\sigma$, then the theta lift of $\sigma$ to $GSO(V)$ is nonzero by the above lemma and all its irreducible constituents are nearly equivalent to $JL(\pi_D) \boxtimes \mu$.
By the generalized strong multiplicity one theorem for $GL_4$, this contradicts the
fact that $JL(\pi_D)$ is non-cuspidal.
\end{proof}

\vskip 15pt

We now consider the local situation.
Over a local field, one has the analogous Weil representation
$\omega_{D,v}$ for $R_D(F_v)$. If $\pi_D\otimes \mu$ and $\sigma$
are irreducible representations of $GSO(V_D)(F_v)$ and $GSp(W)(F_v)$
respectively, then one says that they correspond under theta
correspondence if
\[  \Hom_{R^0_D(F_v)}(\omega_{D,v},  \sigma \boxtimes (\pi_D\boxtimes \mu)) \ne 0. \]
Necessarily, the central characters of $\sigma$ and $\pi_D\otimes
\mu$ are equal. It is perhaps easier to work with the compactly induced Weil
representation
\[  \Omega_{D,v} = ind_{R_D^0}^{GSp(W) \times GSO(V_D)} \omega_{D,v}. \]
It follows from Frobenius reciprocity that $\pi_D\boxtimes \mu$ and
$\sigma$ correspond if and only if
\[  \Hom_{GSp(W) \times GSO(V_D)}(\Omega_{D,v}, \sigma \boxtimes (\pi_D\otimes \mu)) \ne 0. \]
As usual, given $\pi_D\boxtimes \mu$, the maximal $\pi_D\boxtimes
\mu$-isotypic quotient of $\Omega_{D,v}$ is of the form
\[  (\pi_D\boxtimes \mu) \boxtimes \Theta(\pi_D\boxtimes \mu) \]
for some smooth representation $\Theta(\pi_D\boxtimes \mu)$ of $GSp(W)$. One
knows that $\Theta(\pi_D\boxtimes \mu)$ is of finite length and for
lack of a better terminology, we call $\Theta(\pi_D\boxtimes \mu)$
the {\bf big theta lift} of $\pi_D\boxtimes \mu$.
\vskip 10pt

Set
\[  \text{$\theta(\pi_D\boxtimes \mu)=$ the maximal semisimple quotient
of $\Theta(\pi_D\boxtimes \mu)$;} \]
we call it the {\bf small theta lift} of $\pi_D\boxtimes \mu$.  For the case at hand, one knows the following:

\vskip 5pt

\noindent (i) When $\pi_D \boxtimes \mu$ is supercuspidal,  $\Theta(\pi_D\boxtimes \mu)$ is irreducible if it is nonzero,  regardless of the residual characteristic of $F_v$.   This follows from the corresponding result of Kudla in the isometry case, as shown in [GT, Prop. 2.3 and Lemma 3.1].
\vskip 5pt

\noindent (ii) It is a result of Roberts
[Ro2] that if the Howe conjecture for isometry group holds (for
example, away from residual characteristic two), then
$\theta(\pi_D\boxtimes \mu)$ is zero or irreducible. Actually, the
results of Roberts a priori only apply if we are working with representations
of $GO(V_D)$. The passage from $GO(V_D)$ to $GSO(V_D)$ is explained in [GT, Lemma 3.1].
\vskip 5pt

\noindent (iii) In fact, for the representations $\pi_D \boxtimes \mu$ of $GSO(V_D)$ considered in this paper, we show in Section \ref{S:explicit} that $\theta(\pi_D \boxtimes \mu)$ is zero or irreducible, regardless of residual characteristic of $F_v$.
\vskip 10pt

Similarly, starting with the representation $\sigma$ of $GSp_4$,
one has the representations $\Theta_D(\sigma)$ and
$\theta_D(\sigma)$ of $GSO(V_D)$. For the same reasons as above, one knows that the small theta lift $\theta_D(\sigma)$ is either zero or irreducible (see [GT, Prop. 2.3 and Lemma 3.1]).
\vskip 10pt

The above discussion can be summarized in the following
``theta lifting diagram".
\[
\xymatrix{
    GL_2(D) \times GL_1 \ar@{->>}[r]& GSO(V_D)\ar@{<->}[rr]^{JL}\ar@<1ex>[dr]^{\theta}&& GSO(V)\ar@<-1ex>[dl]_{\theta}& GL_4\times GL_1\ar@{->>}[l]\\
    &&GSp_4\ar@<1ex>[ul]^{\theta_D}\ar@<-1ex>[ur]_{\theta_V}&&\\
    D^{\times}\times D^{\times}\ar@{->>}[r]& GSO(D)\ar@{<->}[rr]^{JL}\ar[ur]_{\theta}&&
    GSO(2,2)\ar[ul]^{\theta}&GL_2\times GL_2\ar@{->>}[l]
}
\]
We shall denote all the theta lifts to $GSp_4$ by $\theta$ and the theta
lift from $GSp_4$ to $GSO(V_D)$ and $GSO(V)$ by $\theta_D$ and
$\theta_V$ respectively. Moreover, $JL$ indicates the
Jacquet-Langlands transfer. \\

\vskip 5pt

We conclude this section with a brief discussion on the
functoriality of the above theta correspondence for spherical
representations. The L-group  of $GSp_4$ is $GSp_4(\C)$ and so an
unramified representation of $GSp_4$ corresponds to a semisimple
class in $GSp_4(\C)$. On the other hand, the L-group of $GSO(V)$ is
the subgroup of $GL_4(\C) \times GL_1(\C)$ given by
\[  \{(g,z) \in GL_4(\C) \times GL_1(\C): \det(g) = z^2 \}. \]
There is a natural map
\[  {^L}GSp_4  \longrightarrow  {^L}GSO(V) \]
given by
\[  g \mapsto (g, \Lambda(g)) \]
where $\Lambda$ is the similitude factor of $GSp_4(\C)$.  The following proposition is shown in [GT, Cor. 12.14]:
\vskip 5pt

\begin{Prop}  \label{P:unramGT}
If  $\sigma_v$ is the unramified representation of $GSp_4$ corresponding to the semisimple class $s \in GSp_4(\C)$, then the small theta lift of $\sigma_v$ is the unramified representation of $GSO(V)$ corresponding to the  semisimple class $(s_v, \Lambda(s_v)) \in {^L}GSO(V)$.
\end{Prop}

\vskip 5pt

Moreover, corresponding to the inclusion  $SO(V) \hookrightarrow
GSO(V)$, one has a map of L-groups
\[  std: {^L}GSO(V) \longrightarrow {^L}SO(V) = SO_6(\C). \]
Indeed, one has the map
\[  GL_4(\C) \times GL_1(\C) \longrightarrow GSO_6(\C) \]
given by:
\[  (g,z) \mapsto  z^{-1} \cdot \wedge^2 g, \]
and the map $std$ is simply the restriction of this map to the
subgroup ${^L}GSO(V)$. Thus, one may consider the (partial) standard
degree 6 L-function of a cuspidal representation $\pi \boxtimes \mu$
of $GSO(V)$, which we denote by $L^S(s, \pi \boxtimes \mu, std)$. If
we regard $\pi \boxtimes \mu$ as a representation of $GL_4 \times
GL_1$, then this L-function is nothing but the twisted exterior
square L-function $L^S(s, \pi, \bigwedge^2 \otimes \mu^{-1})$.
\vskip 5pt

Observe finally that if we consider the composite
\[  GSp_4(\C) \longrightarrow {^L}GSO(V) \longrightarrow SO_6(\C), \]
then this 6-dimensional representation of $GSp_4(\C)$ decomposes as
the sum of the trivial representation and the standard 5-dimensional
representation
\[  GSp_4(\C) \longrightarrow PGSp_4(\C) \cong SO_5(\C). \]
Indeed, one has the commutative diagram: \vskip 5pt
\[  \begin{CD}
GSp_4(\C) @>>> {^L}GSO(V) @>>> GL_4(\C) \times GL_1(\C) \\
@VVV   @VVV  @VVV \\
SO_5(\C) @>>> SO_6(\C) @>>> GSO_6(\C)
\end{CD} \]
\vskip 10pt

In view of this,  Prop. \ref{P:unramGT} immediately gives:
\vskip 5pt
\begin{Prop} \label{P:L-identity}
If $\Theta(\pi \boxtimes \mu)$ is cuspidal and contains $\sigma$ as an irreducible constituent, then
\[  L^S(s, \pi, \bigwedge^2 \otimes \mu^{-1}) = L^S(s, \pi \boxtimes \mu, std) = \zeta^S(s) \cdot  L^S(s, \sigma, std), \]
where $L^S(s,\sigma, std)$ is the (partial) standard degree 5
L-function of $\sigma$.
\end{Prop}

\vskip 15pt

\section{\bf The Implication (A) $\Longrightarrow$ (B)}\label{S:a_implies_b}

In this section, we prove the implication (A)$\Longrightarrow$(B) of Thm. \ref{T:main}. In particular, we give a very short proof of the Jacquet-Martin theorem without the 3 conditions present there.
\vskip 10pt

For a nondegenerate character $\chi$ on the unipotent radical $U$ of
a Borel subgroup of $GSp_4$, let $\mathcal{W}_{\chi}$  denote the global
Whittaker functional on $\mathcal{A}(GSp_4)$:
\[  \mathcal{W}_{\chi}: F \mapsto \int_{U(F) \backslash U(\A)} F(u) \cdot \overline{\chi(u)} \, du \]
\vskip 5pt

The following proposition is the key computation (see \cite{So} for
the same computation when $D$ is split):

\vskip 5pt

\begin{Prop}  \label{P:main}
Let $\pi_D$ be a cuspidal representation of $GL_2(D)$ whose central
character is $\mu^2$, so that we may consider the representation
$\pi_D \boxtimes \mu$ of $GSO(V_D)$. Then we have:
\[  \mathcal{W}_{\chi}(\theta(\phi,f)) =  \int_{S_D(\A) \backslash SO(V_D)(\A)} \overline{\mathcal{S}_D( h \cdot f)} \cdot  \left( \int_{U_Y(\A)}
\overline{\chi(u)}  \cdot \omega(u,h)\phi(\underline{x}_0) \, du
\right) \, dh,
\]
where $S_D$ is embedded in to $SO(V_D)$ via $\iota$  and
\[  \underline{x}_0 = (x^0_1, x^0_2) \in V_D(\A)^2 \]
with
\[ \begin{cases}
 x_1^0 = ((1,0);0_D) \in \mathbb{H} \oplus D = V_D \\
x_2^0 = ((0,0); 1_D) \in \mathbb{H} \oplus D = V_D. \end{cases} \]
\end{Prop}
\vskip 5pt

\begin{proof}
Let us write $U = N_Y \rtimes U_Y$ with $U_Y$ a maximal unipotent
subgroup of $GL(Y)$. Then we may restrict the character $\chi$ to
$N_Y$ and $U_Y$. Its restriction to $N_Y$ is a degenerate character,
whereas its restriction to $U_Y$ is nondegenerate. \vskip 5pt

Now we have:
\begin{align}
&\mathcal{W}_{\chi}(\theta(\phi, f))  \notag \\
= &\int_{U(F) \backslash U(\A)} \overline{\chi(u)} \cdot
\int_{SO(V_D)(F) \backslash SO(V_D)(\A)}
\theta(\phi)(u,h)  \cdot \overline{f(h)} \, dh\, du \notag \\
= & \int_{SO(V_D)(F) \backslash SO(V_D)(\A)} \overline{f(h)} \cdot
\int_{U_Y(F) \backslash U_Y(\A)}
 \overline{\chi(u)} \cdot \int_{N_Y(F) \backslash N_Y(\A)}
\overline{\chi(n)} \cdot  \sum_{\underline{x} \in (X \otimes V_D)(F)} \omega(nu,h) \phi (\underline{x})  \, dn \, du \, dh \notag \\
=& \int_{SO(V_D)(F) \backslash SO(V_D)(\A)} \overline{f(h)} \cdot
\left(  \int_{U_Y(F) \backslash U_Y(\A)}
 \overline{\chi(u)} \cdot \sum_{\underline{x} \in \Omega}  \omega(u,h)\phi(\underline{x}) \, du \right) \, dh \notag
\end{align}
where
\[  \Omega = \{ \underline{x} = (x_1,x_2) \in V_D^2: \text{$q_D(x_1) = 0$, $q_D(x_2) = 1$ and $x_2 \in x_1^{\perp}$} \}. \]
Clearly, one has a decomposition
\[  \Omega = \Omega_0 \bigsqcup \Omega_1 \]
where $\Omega_0$ (resp. $\Omega_1$) is the subset of elements with
$x_1 = 0$ (resp. $x_1 \ne 0$). It is easy to see that the sum over
$\Omega_0$ does not contribute and so we need only consider the sum
over $\Omega_1$ above. \vskip 10pt

Now the element $\underline{x}_0$ lies in $\Omega_1$ and the group
$SO(V_D)(F)$ acts transitively on $\Omega_1$.  Indeed, if we set
\[   \Xi = \{  \underline{x}_t = (x_1^0, x_2^0+t x_1^0) : t \in F \} \subset \Omega_1, \]
then
\[   \Omega_1 = SO(V_D)(F) \times_{S_D(F)} \Xi. \]
 Moreover, identifying $\Xi$ with $F$ in the obvious way, the action of $U_Y(F) \times S_D(F)$
 on $\Xi$ is given by:
\[   (u , d \cdot n) : t \mapsto  (u + Tr_D(n)) \cdot t. \]
Hence $U_Y(F)$ acts simply transitively on $\Xi$. \vskip 10pt

Thus we have:
\begin{align}
&\mathcal{W}_{\chi}(\theta(\phi, f))  \notag \\
 =&\int_{S_D(F) \backslash SO(V_D)(\A)} \overline{f(h)} \cdot  \int_{U_Y(F) \backslash U_Y(\A)}
 \overline{\chi(u)} \cdot  \left( \sum_{t \in F}
 \omega(u,h)\phi( \underline{x}_t) \right)  \, du \, dh \notag \\
 =&\int_{S_D(F) \backslash SO(V_D)(\A)}  \overline{f(h)} \cdot \int_{U_Y(\A)}  \overline{\chi(u)} \cdot \omega(u,h)\phi(\underline{x}_0)\, du  \, dh \notag \\
 =& \int_{S_D(\A) \backslash SO(V_D)(\A)} \int_{S_D(F) \backslash
 S_D(\A)} \overline{f(rh)} \cdot
 \int_{U_Y(\A)}  \overline{\chi(u)} \cdot \omega(u, rh)\phi(\underline{x}_0) \, du \, dr \, dh \notag \\
 =& \int_{S_D(\A) \backslash SO(V_D)(\A)} \int_{S_D(F) \backslash
 S_D(\A)} \overline{f(rh)} \cdot
 \int_{U_Y(\A)} \overline{\chi(u)} \cdot \omega(h)\phi(\underline{x}_{u + Tr_D(r)}) \, du \,dr \,  dh \notag \\
 =& \int_{S_D(\A) \backslash SO(V_D)(\A)} \int_{S_D(F) \backslash
 S_D(\A)} \overline{f(rh)}  \chi(Tr_D(r)) \cdot
 \int_{U_Y(\A)}  \overline{\chi(u)} \omega(h)\phi(\underline{x}_{u})  \, du \, dr \, dh  \notag \\
 =& \int_{S_D(\A) \backslash SO(V_D)(\A)} \overline{\mathcal{S}_D( h \cdot f)} \cdot
 \left( \int_{U_Y(\A)}
\overline{\chi(u)}  \cdot \omega(u,h)\phi(\underline{x}_0) \, du
\right) \, dh. \notag
 \end{align}

\end{proof}

\begin{Cor} \label{C:main}
The cuspidal representation $\Theta(\pi_D \boxtimes \mu)$ is
globally generic if and only if $\pi_D$ has Shalika period with
respect to $\mu$.
\end{Cor}
\begin{proof}
This follows from the above proposition by a standard argument
analogous to that in \cite[Pg. 2718-2719]{GS}.
\end{proof}
\vskip 10pt

Now we can prove a part of the implication $(A)\Longrightarrow(B)$, which
is essentially the Jacquet-Martin theorem without the three conditions present there: \vskip 10pt

\begin{Thm}  \label{T:JM2}

If $\pi_D$ has Shalika period with respect to $\mu$, then the
partial L-function $L^S(s, \pi_D, \bigwedge^2 \otimes \mu^{-1})$ has
a pole at $s=1$. Hence its Jacquet-Langlands lift $\pi=JL(\pi_D)$ (if cuspidal)
has Shalika period with respect to $\mu$.
\end{Thm}
\vskip 5pt

\begin{proof}
Suppose that $\pi_D \boxtimes \mu$ has Shalika period. Then by Cor.
\ref{C:main}, $\Theta(\pi_D\boxtimes \mu)$ contains a globally
generic  cuspidal representation $\sigma$ of $GSp_4$. By Prop. \ref{P:L-identity}, we have:
\[  L^S(s,\pi_D, \bigwedge^2 \otimes \mu^{-1})= L^S(s, \pi_D \boxtimes \mu, std) = L^S(s, \sigma, std) \cdot \zeta^S(s). \]
Now  because $\sigma$ is globally generic, $L^S(s,\sigma,
std)$ is non-zero at $s =1$ (by Shahidi \cite[Thm. 5.1]{Sh}). Thus
$L^S(s, \pi_D, \bigwedge^2 \otimes \mu^{-1})$ has a pole at $s =1$.
It follows from results of Jacquet-Shalika \cite{JS} that $\pi$ has
Shalika period with respect to $\mu$.

\end{proof}

\vskip 10pt

Finally the following proposition completes the proof of the implication
$(A)\Longrightarrow(B)$ of Thm. \ref{T:main}.

\begin{Prop} \label{P:bad-type}
Suppose that $D$ is split at every archimedean place.  If $\pi_D$ has
Shalika period with respect to $\mu$, then for all $v \in \Sigma_D$,
the local representation $\pi_{D,v}$ is not equal to
$Ind^{GL_2(D_v)}_{P_{D,v}} \delta_{P_D}^{1/2} \cdot (\tau_{D,1,v} \boxtimes \tau_{D,2,v})$,
   where  $\tau_{D,i,v}$ are representations of $D_v^{\times}$ with central character $\mu_v$.
 \end{Prop}

\begin{proof}
For this, we need one of our local results proved in
Section \ref{S:explicit}. Indeed, if $\pi_D$ has Shalika period with respect to $\mu$,
then  it follows by Cor. \ref{C:main} that  the global theta lift $\theta(\pi_D\boxtimes\mu)$ to
$GSp_4$ is generic. In particular, the local theta lift
$\theta(\pi_{D,v}\boxtimes\mu_v)$ is generic. However, by Thm. \ref{T:local-theta-D}(i), one sees that
 if $\pi_{D,v}$ (for $v\in \Sigma_D$) is of the ``bad"  type described in the proposition  (which is denoted by $PS(\tau_{D,1,v}, \tau_{D,2,v})$ in Section \ref{S:local} and \ref{S:explicit}),  then the local theta lift is non-generic. With this contradiction, the proposition is proved.
 \end{proof}

 \vskip 15pt

\section{\bf A Counterexample to the Converse of Jacquet-Martin}\label{S:counter_ex}

Before proving the other implication $(B)\Longrightarrow (A)$ of our
main theorem, we describe in this section a concrete counterexample to the
converse of the Jacquet-Martin theorem. Namely we shall construct a cuspidal representation
$\pi$ of $PGL_4$ which has cuspidal Jacquet-Langlands lift
$\pi_D$ on $PGL_2(D)$ and show that $\pi$ has Shalika period but
$\pi_D$ does not. To do so, we shall construct an irreducible cuspidal
representation $\sigma$ on $PGSp_4$ with the following properties:
\vskip 10pt

\begin{itemize}

\item[(i)] $\sigma$ is globally generic.
\vskip 5pt

\item[(ii)] at two finite places $v_1$ and $v_2$, $\sigma_{v_i}$ is supercuspidal and
has a non-zero theta lift to the split $PGSO(2,2)$, say  $\sigma_{v_i} = \Theta(\tau_i\boxtimes \tau'_i)$
for supercuspidal
$\tau_i\boxtimes \tau'_i$ on $PGSO(2,2) = PGL_2 \times PGL_2$.
\vskip 5pt

\item[(iii)] the global theta lift of $\sigma$ to $GSO(2,2)(\A)$ is zero.
\end{itemize}

 \vskip 10pt

Before showing how to construct such a $\sigma$,  let us see why
having such a $\sigma$ gives a counterexample. Because $\sigma$ is
globally generic and cuspidal,  $\pi := \Theta_V(\sigma)$ is nonzero on $PGL_4$ by Lemma \ref{L:generic_non-vanish}. Moreover,  $\pi$ is cuspidal: this follows from
the tower property of theta lifts, since the theta lift of $\sigma$
to $PGSO(2,2)$ is zero by assumption. By the strong multiplicity one
theorem, we see that $\pi$ is irreducible. For $i =1$ or $2$,  a
simple calculation of local theta correspondence (cf. Thm. \ref{T:local-theta}(i)) shows that
\[  \pi_{v_i} =  PS(\tau_i, \tau'_i) := Ind_P^{GL_4} \delta_P^{1/2} \cdot (\tau_i \boxtimes  \tau'_i). \]
\vskip 5pt \noindent Moreover, we have
\[
    L^S(s,\pi, \bigwedge^2 ) = L^S(s, \sigma, std) \cdot \zeta^S(s),
\]
and $L^S(s, \sigma, std)$ is nonzero at $s=1$ by genericity of $\sigma$. Thus $L^S(s,\pi, \bigwedge^2 )$ has a pole at $s=1$, and so $\pi$
has Shalika period. \vskip 5pt

Now let D be the quaternion algebra ramified precisely at $v_1$ and
$v_2$. By recent results of Badalescu \cite{B}, $\pi$ has a cuspidal Jacquet-Langlands lift $\pi_D$ on $PGL_2(D)$ (since $\pi_{v_i}$ is compatible with $GL_2(D_{v_i})$ for $i = 1$ and $2$ in the sense of [B]). We need to show that $\pi_D$ has no Shalika period.
\vskip 5pt

For $i=1$ and $2$,  we have:
\[  \pi_{D,v_i} = PS(JL(\tau_i), JL(\tau'_i)): = Ind_{P_D}^{GL_2(D)}  JL(\tau_i) \boxtimes JL(\tau'_i). \]
In \cite[Prop. 7]{P1}, Prasad showed that
\[    (PS(JL(\tau_i), JL(\tau'_i)))_{N_D, \psi_D} \cong JL(\tau_i) \boxtimes JL(\tau'_i) \]
as representations of $PD_{v_i}^{\times}$, and $\pi_{D,v_i}$ does not posses Shalika period if $\tau_i \ne
\tau'_i$; we recall his results in
Prop. \ref{P:prasad} below. Thus, if we assume that $\tau_i \ne
\tau'_i$ for $i=1$ and $2$, then the local components of $\pi_D$ at
$v_1$ and $v_2$ do not possess local Shalika periods.  Hence $\pi_D$
does not possess Shalika period.

\vskip 5pt

In fact, we can give another argument to show that the same
conclusion holds even if $\tau_i = \tau'_i$ (in which case
$\pi_{D,v_i}$ does possess local Shalika period). Suppose for the
sake of contradiction that $\pi_D$ has Shalika period. Then as in
the proof of Thm. \ref{T:JM2}, its theta lift  $\Theta_D(\pi_D)$ to
$G = PGSp_4$ is non-zero cuspidal and globally generic. Moreover,
any of its irreducible constituents is nearly equivalent to $\sigma$. By the strong multiplicity one
theorem for globally generic cuspidal representations of $PGSp_4$
\cite[Cor. 2]{J-So}, we have
\[  \sigma \subset \Theta(\pi_D).  \]
\vskip 5pt

\noindent This implies that for all $v$, $\sigma_v$  has
non-zero theta lift to $PGSO(V_D)$. On the other hand, Kudla-Rallis
\cite[Thm. 3.9]{KR1} showed that if a supercuspidal representation
of $PGSp_4$ lifts to $PGSO(2,2)$, then it does not lift to
$PGSO(V_D)$. So at the two places $v_1$ and $v_2$, the theta lift of
$\sigma_{v_i}$ to $PGSO(V_D)$ is zero. This gives the desired
contradiction. \vskip 15pt

It remains then to construct a cuspidal representation $\sigma$ with
the properties (i)-(iii) above. Such a $\sigma$ can be obtained as a
theta lift from a cuspidal representation of a quasi-split
$PGSO(3,1)$. Choose a quadratic extension $E/F$ such that $v_1$ and
$v_2$ split in $E$. Let $V_E$ be the quadratic space $(E
,\mathbb{N}_E) \oplus \mathbb{H}$.  Then (cf. \cite[\S 2]{Ro1})
\[  PGSO(V_E)  \cong PGL_2(E).  \]
Let $\Sigma$ be a cuspidal representation of  $PGL_2(E)$ satisfying
the following properties: \vskip 5pt

\noindent (a) for $i = 1$ or $2$, $\Sigma_{v_i}$ is a  supercuspidal
representation of $PGL_2(E_{v_i}) = PGL_2(F_{v_i}) \times
PGL_2(F_{v_i})$. \vskip 5pt

\noindent (b) at some finite place $v_3$ which is inert in $E$, $\Sigma_{v_3}$ is a
supercuspidal representation which is not the base change lift of a
representation of $GL_2(F_{v_3})$.  \vskip 5pt

\noindent Such a $\Sigma$ exists by Poincare series techniques, as shown in [Sh2]. \vskip 5pt

Now consider the theta lift of $\Sigma$ from $PGSO(V_E)$ to
$PGSp_4$, which has been studied in great detail by B. Roberts
\cite{Ro1}. One knows that the theta lift of $\Sigma$ to $PGSp_4$ is
non-zero cuspidal (because of condition (b)) and contains a globally
generic constituent. We take such a constituent  to be our $\sigma$.
By construction, $\sigma$ satisfies the requirements (i) and (ii).
To show that the global theta lift of $\sigma$ to $GSO(2,2)$
vanishes, one considers the standard degree 5 L-function of
$\sigma$. Since $\sigma$ is a theta lift from $PGSO(V_E)$, we have
\[  L^S(s,\sigma, std) = L^S(s, \omega_{E/F}) \cdot L^S(s, \Sigma, Asai) \]
where $\omega_{E/F}$ is the quadratic character associated to $E/F$
and the last L-function is the Asai L-function of $\Sigma$. It is
well-known that $L^S(s, \omega_{E/F})$ is nonzero holomorphic at
$s=1$. Also one knows that $L^S(s, \Sigma, Asai)$ is nonzero
holomorphic at $s=1$ (see \cite{F}, especially \S 5.) Hence the RHS
is nonzero holomorphic at $s =1$. On the other hand, if the theta
lift of $\sigma$ to $GSO(2,2)$ is nonzero, say $\Theta(\sigma)
\supset \tau_1 \boxtimes \tau_2$, then
\[  L^S(s,\sigma, std) = \zeta^S(s) \cdot L^S(s, \tau_1 \times \tau_2), \]
which has a pole at $s =1$. This contradiction shows that $\sigma$
satisfies the condition (iii).
 \vskip 15pt

This completes the construction of the counterexample. Indeed, by
our construction above, we obtain: \vskip 10pt

\begin{Thm}
Suppose that $D$ is split at every archimedean place. Then there are
cuspidal representations $\pi_D$ of $PGL_2(D)$, with a cuspidal
Jacquet-Langlands lift to $PGL_4$, satisfying: \vskip 5pt
\begin{itemize}
\item[(i)] for each place $v$, $\pi_{D,v}$ has local Shalika period, and
\vskip 5pt

\item[(ii)] $L^S(s, \pi_D, \bigwedge^2)$ has a pole at $s= 1$, but
\vskip 5pt

\item[(iii)]  $\pi_D$ does not have global Shalika period.
\end{itemize}
\end{Thm}

\vskip 5pt

\begin{proof}
As in the construction of the counterexample above, we pick a
quadratic field extension $E/F$ such that for all $v \in \Sigma_D$,
$E_v$ is split. For all $v \in \Sigma_D$, pick supercuspidal
representations $\Sigma_v = \tau_v \boxtimes \tau_v$ of $PGL_2(F_v)
\times PGL_2(F_v)$. At some finite place $w \notin  \Sigma_D$ such
that $E_w$ is a field, pick a supercuspidal representation
$\Sigma_w$ of $PGL_2(E_w)$ which is not in the image of the base
change from $GL_2(F_w)$. Let $\Sigma$ be a cuspidal representation
of $PGSO(V_E)$ with these local components and let $\sigma$ be a
globally generic constituent of the theta lift of $\Sigma$ to
$PGSp_4$. As above, if $\pi$ is the global theta lift of $\sigma$ to
$PGSO(V) \cong PGL_4$, then $\pi$ is a cuspidal representation with
Shalika period whose local components at $v \in \Sigma_D$ are of the
form $PS(\tau_v, \tau_v)$. In particular, $L^S(s, \pi, \bigwedge^2)$
has a pole at $s=1$. By Badulescu \cite{B}, this $\pi$ has a
Jacquet-Langlands transfer $\pi_D$ to $PGL_2(D)$ and for $v \in
\Sigma_D$,
\[
    \pi_{D,v} = PS(JL(\tau_v), JL(\tau_v)),
 \]
which has local Shalika period by Prop. \ref{P:prasad}. Then, as we have seen above, all the local components of $\pi_D$ have
Shalika period, but $\pi_D$ itself does not.
\end{proof}
\vskip 15pt

\section{\bf The Implication (B) $\Longrightarrow$ (A)}\label{S:b_implies_a}

In this section we will prove the other implication $(B)\Longrightarrow
(A)$ in Thm. \ref{T:main}. For this, we shall need some precise
results about the global Jacquet-Langlands correspondence between
$GL_2(D)$ and $GL_4$. Such results have now been obtained by
Badulescu \cite{B} in essentially complete generality. However, for technical
reasons, Badulescu \cite{B} assumes that the quaternion algebra $D$
is split at every archimedean place. Thus, at times, we shall need
to make this assumption in this section.

\vskip 10pt

We first note the following lemma: \vskip 5pt

\begin{Lem}
Assume that $D$ is split at every archimedean place. If the cuspidal
representation $\pi =JL(\pi_D)$ of $GL_4$ has Shalika model with
respect to $\mu$, so that $\Theta(\pi \boxtimes\mu)$ contains an irreducible
globally generic cuspidal representation $\sigma$ of $GSp_4$, then the following are
equivalent: \vskip 5pt

(i)  $\pi_D$ has Shalika period with respect to $\mu$; \vskip 5pt

(ii)  the theta lift of $\sigma$ to $GSO(V_D)$ is non-zero. \vskip 5pt

\noindent When these conditions hold, the theta lift of $\sigma$ to $GSO(V_D)$
is equal to $\pi_D \boxtimes \mu$.
\end{Lem}

\begin{proof}
We first prove that (i) implies (ii).
If $\pi_D$ has Shalika period with respect to $\mu$, then by
Corollary \ref{C:main}, its theta lift to $GSp_4$ contains an irreducible
cuspidal globally generic representation $\sigma'$.
Moreover, $\sigma'$ and $\sigma$ are
nearly equivalent and so are equal by the results of Jiang-Soudry
\cite{J-So}. The theta lift $\Theta_D(\sigma)$ of $\sigma$ to $GSO(V_D)$ is
thus nonzero cuspidal and all its constituents are nearly equivalent to $\pi_D \boxtimes \mu$;
it  is thus equal to $\pi_D\boxtimes \mu$ by
the strong multiplicity one theorem for $GL_2(D)$ due to Badulescu
\cite[Thm. 5.1 (b) and (c)]{B}.
\vskip 10pt

Conversely, if the theta lift of $\sigma$ to $GSO(V_D)$ is non-zero,
then all its irreducible constituents are nearly equivalent to $\pi_D\boxtimes \mu$ and thus $\Theta_D(\sigma)$ is
equal to $\pi_D \boxtimes \mu$ by the strong multiplicity one
theorem for $GL_2(D)$ \cite[Thm. 5.1 (b)
and (c)]{B}. This shows that the theta lift of $\pi_D \boxtimes \mu$ to $GSp_4$ is not orthogonal to $\sigma$, and thus contains an irreducible constituent $\sigma'$ isomorphic to $\sigma$. Again, the results of \cite{J-So} implies that $\sigma' = \sigma$, so that $\Theta(\pi_D \boxtimes \mu)$ is globally generic. Corollary \ref{C:main} then implies that $\pi_D \boxtimes\mu$ has Shalika period with respect to $\mu$.
\end{proof}

\vskip 5pt

Thus, a necessary condition for $\pi_D$ to have Shalika period is
that the local representations $\sigma_v$ has a non-zero theta lift
to $GSO(V_D)(F_v)$. This is of course automatic for $v \notin
\Sigma_D$, but is not automatic for $v \in \Sigma_D$ (as the counterexample shows),
and hence we need the local condition as in our main
theorem. Of course, even when these local obstructions to theta
lifting are absent, one still has to show that the global theta lift
is non-zero.

\vskip 10pt In any case, the following theorem immediately implies
the implication $(B)\Longrightarrow (A)$: \vskip 5pt

\begin{Thm} \label{T:converse}
Suppose that $\pi_D$ is a cuspidal representation of $GL_2(D)$ with
central character $\mu^2$ and a cuspidal Jacquet-Langlands lift
$\pi$ on $GL_4$.  If \vskip 5pt

\noindent (i) the $L$-function $L^S(s, \pi_D,
\bigwedge^2\otimes\mu^{-1})$ has a pole at $s=1$, \ie $\pi$ has
Shalika period with respect to $\mu$, and \vskip 5pt

\noindent (ii)  for $v \in \Sigma_D$, $\pi_v$ is not of the form
$PS(\tau_{1,v} \boxtimes \tau_{2,v})$ where $\tau_{i,v}$ are
representations of $GL_2(F_v)$ with central character $\mu$, \vskip
5pt

\noindent then there is a cuspidal representation $\pi'_D$ on
$GL_2(D)$ which is nearly equivalent to $\pi_D$ and which possesses
Shalika period with respect to $\mu$. Moreover, if $D$ is split at
every archimedean place of $F$, then $\pi'_D$ is equal to $\pi_D$.
\end{Thm}

\begin{proof}
Let $\sigma = \Theta(\pi \boxtimes \mu)$. By (i), $\sigma$ is
globally generic and irreducible cuspidal. By (ii), for $v \in \Sigma_D$, the
local components $\sigma_v$ do not participate in the local theta
correspondence with $GSO(2,2)$.  In view of the previous lemma, in
order to prove the theorem, it suffices to show that $\pi'_D
:=\Theta_D(\sigma) \ne 0$. As we explain in Section \ref{S:theta}, this non-vanishing is equivalent to
the non-vanishing of the global theta lift of $\sigma|_{Sp_4}$ to $O(V_D)$ (where we are considering the restriction of functions from $GSp_4$ to $Sp_4$). Thus, we may work with isometry groups below.
We shall in fact show that the period of
$\Theta_D(\sigma|_{Sp_4})$ over the subgroup $O(D) \subset O(V_D)$ is non-zero. Our
argument below is largely inspired by \cite[\S 7]{KRS}. \vskip 10pt

The decomposition $V_D = D \oplus \mathbb{H}$ gives a see-saw
diagram: \vskip 5pt

\[
\xymatrix{ Sp_4\times Sp_4\ar@{-}[dr] & O(V_D)\ar@{-}[dl]\\
\Delta Sp_4 \ar@{-}[u]  & O(D)\times O(\mathbb{H})\ar@{-}[u] }
\]

\vskip 5pt On restriction to $\Delta Sp_4 \times (O(D) \times
O(\mathbb{H}))$, the Weil representation of $Sp_4 \times O(V_D)$
decomposes as a tensor product:
\[  S(X \otimes V_D) \cong S(X \otimes D) \hat{\otimes} S(X \otimes \mathbb{H}). \]
Let $\phi = \phi_D \otimes \phi_{\mathbb{H}} \in S((X \otimes
V_D)(\A))$ and take $f\in\sigma$. Now let us  compute the period
of $\theta_D(\phi,f)$ over the anisotropic  group $O(D)$. We get:
\vskip 5pt
\begin{align}
&\int_{O(D)(F)  \backslash O(D)(\A)} \theta_D(\phi, f)(h) \, dh  \notag \\
= & \int_{O(D)(F) \backslash O(D)(\A)}  \int_{Sp_4(F) \backslash Sp_4(\A)} \theta(\phi)(gh) \cdot \overline{f(g)} \, dg \, dh \notag \\
=& \int_{O(D)(F) \backslash O(D)(\A)}  \int_{Sp_4(F) \backslash Sp_4(\A)} \theta(\phi_D)(gh) \cdot \theta(\phi_{\mathbb{H}})(g) \cdot \overline{f(g)} \, dg \, dh \notag \\
=& \int_{Sp_4(F) \backslash Sp_4(\A)} \overline{f(g)}  \cdot
\theta(\phi_{\mathbb{H}})(g) \cdot \left( \int_{O(D)(F) \backslash
O(D)(\A)} \theta(\phi_D)(gh) \, dh \right) \,  dg. \notag
\end{align}
\vskip 5pt \noindent The inner integral is the theta lift of the
trivial representation of $O(D)$ to $Sp_4$. Note that since $O(D)$ is anisotoropic, this theta integral is always convergent and hence there is no need for regularization. Now by the Siegel-Weil
formula of Kudla-Rallis \cite{KR2}, this inner integral is equal to
an Eisenstein series described as follows. There is a
$O(D)$-invariant and $Sp_4$-equivariant map
\[  F: S((X \otimes D)(\A)) \longrightarrow I(1/2) \]
with
\[  I(s) = Ind_{P(Y)}^{Sp_4} |\det|^s \quad \text{(normalized induction)} \]
given by
\[  F(\phi_D)(g) = \omega(g)\phi_D(0). \]
One may consider the Eisenstein series $E(F(\phi_D),s,g)$ associated
to the standard section attached to $F(\phi_D)$.  Then the result of
Kudla-Rallis is:
\[    \int_{O(D)(F) \backslash O(D)(\A)} \theta(\phi_D)(gh) \, dh  = c \cdot E(F(\phi_D), 1/2, g) \]
for some non-zero constant $c$. By adjusting the measure $dh$, there
is no loss of generality in assuming that $c =1$. \vskip 10pt

It should be noted that the family of Eisenstein series attached to
an arbitrary  standard section of $I(s)$ can have a pole of order
$1$ at $s = 1/2$. However, for the sections in the image of $F$, the
associated Eisenstein series is holomorphic at $s = 1/2$. This is
reflected by the fact that $I(1/2)$ is reducible. The structure of
the local degenerate principal series $I_v(1/2)$ is described
precisely in \cite[Props. 1.1 and 1.2]{KRS} and \cite[Thm. 1]{LZ}. We record the relevant
facts: \vskip 5pt

\begin{Prop}
(i) If $v$ is non-archimedean, then $I_v(1/2) =
\Theta_v(1_{O(2,2)})$, which has length 2. It has a unique
irreducible submodule isomorphic to $\Theta_v(1_{O(D)})$ and a
unique irreducible quotient isomorphic to $\Theta_v(1_{O(1,1)})$.
\vskip 5pt

(ii) If $v$ is archimedean, then $I_v(1/2) = \Theta_v(1_{O(2,2)})$.
If $v$ is real, then $I_v(1/2)$ has a unique irreducible quotient isomorphic to
$\Theta_v(1_{O(1,1)})$ and its unique maximal submodule is isomorphic to
$\Theta_v(1_{O(D)}) := \Theta_v(1_{O(4,0)}) \oplus \Theta_v(1_{O(0,4)})$.
\end{Prop}
\vskip 5pt

\begin{Cor}
The image of $F$ is the submodule of $I(1/2)$ given by
\[  \left( \otimes_{v \in \Sigma_D} \Theta_v(1_{O(D)})  \right) \bigotimes
\left( \otimes_{v \notin \Sigma_D} I_v(1/2). \right) \]
\end{Cor}
\vskip 10pt

In view of the Siegel-Weil formula, we see that
\[ \int_{O(D)(F)  \backslash O(D)(\A)} \theta_D(\phi, f)(h) \, dh  =
 \int_{Sp_4(F) \backslash Sp_4(\A)} \overline{f(g)}  \cdot \theta(\phi_{\mathbb{H}})(g) \cdot
E(F(\phi_D), 1/2,g) \,  dg \] and to prove the theorem, it suffices
to show that the integral on the RHS is non-zero for some choices of
$\phi= \phi_D \otimes \phi_{\mathbb{H}}$ and $f$. \vskip 10pt

In \cite{PSR}, Piatetski-Shapiro and Rallis have considered the
Rankin-Selberg integral suggested by the RHS of the above equality:
\[  Z(s, f,  \Phi, \phi_{\mathbb{H}}) =
 \int_{Sp_4(F) \backslash Sp_4(\A)} \overline{f(g)}  \cdot \theta(\phi_{\mathbb{H}})(g) \cdot
E(\Phi, 1/2,g) \,  dg. \] 
This family of global zeta integrals is
not identically zero if $\sigma$ has non-vanishing Fourier
coefficients along $N_Y$  corresponding to the split binary quadratic
space. This is the case since $\sigma$ is globally generic.
Piatetski-Shapiro and Rallis showed that
\[  Z(s,f,\Phi, \phi_{\mathbb{H}}) = L^S(s+ \frac{1}{2}, \sigma, std) \cdot  Z_S(s,f _S,\Phi_S,\phi_{\mathbb{H},S})\]
where $S$ is a finite set of places of $F$ containing $\Sigma_D$.
\vskip 5pt

Now let us examine the analytic behavior of both sides at $s=1/2$.
The Eisenstein series $E(\Phi,s,g)$ has a pole of order at most 1 at
$s=1/2$ and its residue there is contained in the regularized theta lift of
the trivial representation of $O(1,1)(\A)$ \cite[Thm.
4.1(iii)]{KRS}. Thus, if the residue at $s=1/2$ of the LHS is
nonzero, we would conclude that $\sigma$ has a non-zero theta lift
to $GO(2,2)$, which is a contradiction (since we know that $\sigma_v$
does not lift to $GO(2,2)$ for $v \in \Sigma_D$). Thus, the LHS is
holomorphic at $s=1/2$. On the other hand, we know that
$L^S(s+\frac{1}{2},\sigma, std)$  is holomorphic and nonzero at
$s=1/2$. This implies that the ramified factor $Z_S(s,f_S,\Phi_S,
\phi_{\mathbb{H},S})$ is also holomorphic at $s =1/2$. Moreover, it
was shown in \cite{PSR} that there are choices of data such that
$Z_S(1/2, f_S,\Phi_S,\phi_{\mathbb{H},S})$ is nonzero. 
\vskip 10pt

For our purpose, we need to show that for some $\Phi$ of the form
$F(\phi_D)$, the ramified factor $Z_S(1/2, f_S,\Phi_S,
\phi_{\mathbb{H},S}) \ne 0$. Let us first fix $f_S$,  $\Phi_S$ and
$\phi_{\mathbb{H},S}$ such that
$Z_S(1/2,f_S,\Phi_S,\phi_{\mathbb{H},S}) \ne 0$. Now, fixing the
components of $f_S$ in $S \setminus \Sigma_D$ while varying the
components in $\Sigma_D$, we see that this ramified zeta factor at $s = 1/2$
gives  a nonzero $Sp_4(F_{\Sigma_D})$-equivariant map
\[   \sigma_{\Sigma_D}^{\vee} \otimes I_{\Sigma_D}(1/2) \otimes S((X \otimes \mathbb{H})(F_{\Sigma_D})) \longrightarrow
\C. \] 
We need to show that it is still non-zero if we restrict the
second argument  to the submodule
$\Theta_{\Sigma_D}(1_{O(D)})$. If not, then for some place $v \in
\Sigma_D$, we would obtain a non-zero $Sp_4(F_v)$-equivariant map
\[  \sigma_v^{\vee} \otimes \Theta_v(1_{O(1,1)}) \otimes S((X \otimes \mathbb{H})(F_v)) \longrightarrow \C. \]
 Since $\Theta_v(1_{O(1,1)})$ is a quotient of $S(X \otimes \mathbb{H})$, we would deduce that
 there is a non-zero $Sp_4(F_v)$-equivariant map
 \[  S((X \otimes \mathbb{H}^2)(F_v)) \longrightarrow \sigma_v. \]
 This contradicts the assumption that $\sigma_v$ does not participate in the theta correspondence with $GO(2,2)$.
 \vskip 10pt

This completes the proof of Theorem \ref{T:converse}.
\end{proof}

\vskip 15pt

\noindent{\bf Remarks:}  Indeed, what the proof of Theorem
\ref{T:converse} shows is the following. Suppose that $\sigma$ is a
cuspidal (not necessarily generic) representation of $GSp_4(\A)$
satisfying: \vskip 5pt

\begin{itemize}
\item $\sigma$  has nonzero Fourier coefficient along $N_Y$ corresponding to the split binary quadratic space;

\item $L^S(1, \sigma , std)$ is finite but nonzero;

\item for all $v \in \Sigma_D$, $\sigma_v$ does not participate in the theta correspondence with $GO(2,2)$.
\end{itemize}

\noindent Then the global theta lift of $\sigma$ to $GSO(V_D)$ is
nonzero.

\vskip 15pt

\section{\bf The Equivalences of (C), (D) and (E)}  \label{S:cde}

In this section, we show the equivalences of (C), (D) and (E) in Thm. \ref{T:main}.
For convenience, we restate the result to be proved:
\vskip 5pt

\begin{Thm}
Assume that $D$ is split at every archimedean place.
Suppose that $\pi_D \boxtimes \mu$ is a cuspidal representation of $GSO(V_D)$
whose Jacquet-Langlands lift $JL(\pi_D)$ is not cuspidal, so that
 $JL(\pi_D)$ is contained in the residual spectrum and is isomorphic
 to the unique irreducible quotient of  $PS(\tau |-|^{1/2}, \tau |-|^{-1/2})$
for a cuspidal representation $\tau$ of $GL_2(\A)$. Then the following are equivalent:
\vskip 5pt
\begin{enumerate}
\item[(C)]  $\pi_D$ has Shalika period with respect to $\mu$.
\vskip 5pt

\item[(D)] $\mu$ is equal to the central character $\omega_{\tau}$ of $\tau$.
\vskip 5pt

\item[(E)]  The (incomplete) twisted
exterior square $L$-function $L^S(s, \pi_D, \bigwedge^2\otimes\mu^{-1})$ has a pole at $s= 2$.
\end{enumerate}
\end{Thm}

\vskip 10pt

The equivalence of (D) and (E) is easy to verify. Indeed, since $\pi_D$ is nearly equivalent to
any irreducible constituent of
$PS(\tau|-|^{1/2}, \tau|-|^{-1/2})$, we see that $\mu^2 = \omega_{\tau}^2$
so that $\mu= \omega_{\tau} \cdot \chi$ for some quadratic character $\chi$ and
\[  L^S(s, \pi_D, \bigwedge^2\otimes\mu^{-1}) = L^S(s+1,\chi) \cdot L^S(s-1, \chi) \cdot L^S(s, \tau \times \tau^{\vee} \cdot \chi^{-1}). \]
From this, one deduces that $L^S(s, \pi_D, \bigwedge^2\otimes\mu^{-1})$ has a pole at $s = 2$ if and only if $\chi$ is trivial, i.e. that $\mu = \omega_{\tau}$.
\vskip 10pt

Before proving the equivalence of (C) and (D), let us take note of the following consequence of Badulescu's paper [B].
\vskip 5pt

\begin{Prop}  \label{P:residue}
The Jacquet-Langlands correspondence (as defined by Badulescu [B]) sets up a bijection between
\vskip 5pt
\begin{itemize}
\item[(a)]  the set of irreducible infinite dimensional constituents of the discrete spectrum of $GL_2(D)$ whose Jacquet-Langlands  lift to $GL_4$ is non-cuspidal;
\vskip 5pt

\item[(b)] the irreducible infinite dimensional constituents of the residual spectrum of $GL_4$.
\end{itemize}
\vskip 5pt

\noindent  Moreover, for  $\pi_D$ as in (a), so that $JL(\pi_D)$ is the unique irreducible quotient of
$PS(\tau|-|^{1/2}, \tau|-|^{-1/2})$ for a cuspidal representation $\tau$ of $GL_2(\A)$, $\pi_D$ is
cuspidal if and only if  $\tau$ is not compatible with $D$,
i.e. $\tau_v$ is not a discrete series representation for some $v \in \Sigma_D$.
\end{Prop}
 \vskip 10pt

Now suppose that (C) holds so that $\pi_D$ has Shalika period with respect to $\mu$. Then
by Cor. \ref{C:main} and Prop. \ref{P:JL-theta}(ii), the theta lift of $\pi_D \boxtimes \mu$ to $GSp_4$ is nonzero and non-cuspidal. Thus, by the tower property of theta correspondence, the theta lift of $\pi_D \boxtimes \mu$ to $GSp_2 \cong GL_2$ is nonzero and cuspidal. Let $\sigma$ be an irreducible constituent of the theta lift of $\pi_D \boxtimes \mu$ to $GL_2$, so that $\omega_{\sigma} = \mu$. Since the theta lift of $\sigma$ to
$GSO(2,2)$ is nonzero, the theta lift of $\sigma$ to $GSO(V)$ is also nonzero. Indeed, it is not difficult to check that the theta lift of $\sigma$ to $GSO(V)$ is nearly equivalent to the irreducible
constituents of $PS(\sigma|-|^{1/2}, \sigma|-|^{-1/2}) \boxtimes \omega_{\sigma}$ (cf. Thm.
\ref{T:explicit}). However, this is nearly equivalent to $JL(\pi_D)$ and so we deduce by the generalized strong multiplicity one theorem
that
\[  \sigma = \tau \quad \text{and} \quad \omega_{\tau} = \mu. \]
This proves the implication (C)$\Longrightarrow$(D).

\vskip 10pt

Suppose now that  (D) holds. Then we consider the theta lift $\Theta_D(\tau)$ of $\tau$ from $GL_2$ to $GSO(V_D)$.
In the proposition below, we shall show that $\Theta_D(\tau)$ has nonzero Shalika period. This is sufficient to show the implication (D)$\Longrightarrow$(C). Indeed, by Prop. \ref{P:residue}, the fact that $\pi_D$ is cuspidal implies that $\tau$ is not compatible with $D$, so that the theta lift of $\tau$ to $GSO(D)$ is zero. This shows that $\Theta_D(\tau)$ is cuspidal and  all its constituents are nearly equivalent to $PS(\tau|-|^{1/2}, \tau|-|^{-1/2}) \boxtimes \omega_{\tau}$ and thus to
$\pi_D \boxtimes \omega_{\tau}$. By the strong multiplicity one theorem for $GL_2(D)$,
one concludes that $\pi_D \boxtimes \omega_{\tau}= \Theta_D(\tau)$ and so $\pi_D$ has Shalika period with respect to $\mu$. This proves (C).
\vskip 10pt

It remains then to show:
\vskip 5pt

\begin{Prop} \label{P:shal}
Let $\tau$ be a cuspidal representation of $GL_2$ and let $\Theta_D(\tau)$ denote the theta lift of $\tau$ to $GSO(V_D)$. Then $\Theta_D(\tau)$ has nonzero Shalika period.
\end{Prop}

\begin{proof}
This follows by a direct computation. We begin by setting up some notations. Let
$W' = F \cdot e \oplus F \cdot f$ be a rank 2 symplectic space so that $GSp(W') \cong GL_2$.
Then $Sp(W')$ acts transitively on the nonzero elements of $W$ and the stabilizer of $e$ is the
unipotent radical $U$ of the Borel subgroup stabilizing the line $F \cdot e$.
\vskip 5pt

Recall that  we have a decomposition
\[  V_D = F \cdot (1,0) \oplus D \oplus F \cdot (0,1) \]
and let us set $v_0 = (1,0)$ and $v_0^* = (0,1)$.
The Weil representation $\omega'_D$ of $O(V_D) \times Sp(W')$ has a mixed model relative to the above decomposition of $V_D$. This is realized on
the space of Schwarz functions  on $(v_0^* \otimes W) \oplus (f \otimes D)$, i.e. on
$S(v_0^* \otimes W) \otimes S(f \otimes D)$.
\vskip 5pt

For  $\phi \in S(v_0^* \otimes W) \otimes S(f \otimes D)$, let $\theta(\phi)$
be the associated theta function and for $f \in \tau$, we have the theta lift  $\theta(\phi, f)$.
 Now we compute:

\begin{align}
&\mathcal{S}_D(\theta(\phi,f))  \notag \\
=&\int_{S_D(F) \backslash S_D(\A)} \overline{\psi_D(s)} \cdot \left(
\int_{Sp(W')(F) \backslash Sp(W')(\A)} \theta(\phi)(sg) \cdot \overline{f(g)} \, dg  \right) \, ds \notag \\
=& \int_{S_D(F) \backslash S_D(\A)} \overline{\psi_D(s)} \cdot  \left(
\int_{Sp(W')(F) \backslash Sp(W')(\A)}  \sum_{w \in W} \, \sum_{d \in D}
(\omega_D'(sg)\phi)(w,  d)  \cdot \overline{f(g)} \, dg \right) \, ds  \notag \\
=&  \int_{S_D(F) \backslash S_D(\A)} \overline{\psi_D(s)} \cdot  \left(
\int_{Sp(W')(F) \backslash Sp(W')(\A)}  \sum_{\gamma \in U(F) \backslash Sp(W')(F)}
\sum_{d \in D} (\omega_D' (s \gamma g)\phi) (e, d) \cdot \overline{f(\gamma g)}  \, dg \right) \, ds \notag \\
=&  \int_{S_D(F) \backslash S_D(\A)} \overline{\psi_D(s)} \cdot
\left( \int_{U(F) \backslash Sp(W')(\A)}
\sum_{d \in D} (\omega_D' (sg)\phi) (e, d) \cdot \overline{f(g)} \, dg  \right) \, ds \notag \\
=& \int_{U(F) \backslash Sp(W')(\A)}  \overline{f(g)} \cdot   \int_{PD^{\times}_F
\backslash PD^{\times}_{\A}}
 \sum_{d \in D}    \omega_D' (hg)\phi (e, d) \cdot
\left( \int_{N_D(F)\backslash N_D(\A)} \overline{\psi_D(n)}  \cdot  \psi(Tr(n \overline{d}))  \, dn \right) \,
dh \, dg  \notag \\
=& \int_{U(F) \backslash Sp(W')(\A)}   \overline{f(g)} \cdot
\int_{PD^{\times}_F \backslash PD^{\times}_{\A}}
   \omega_D' (hg)\phi (e, 1_D)  \, dh \, dg  \notag \\
=& \int_{U(\A) \backslash Sp(W')(\A)}  \left( \int_{U(F) \backslash U(\A)}  \psi(u) \cdot \overline{f(ug)}
 \, du \right)  \cdot
\omega_D'(g)\phi(e, 1_D)   \, dg \notag \\
=& \int_{U(\A) \backslash Sp(W')(\A)}   \overline{\mathcal{W}_{\psi}(g \cdot f)} \cdot
\omega_D'(g)\phi(e, 1_D)   \, dg, \notag
\end{align}
where $\mathcal{W}_{\psi}$ is the global Whittaker functional on $\tau$ and we have normalized measures so that
\[  \int_{N_D(F) \backslash N_D(\A)} dn = 1 \quad \text{and} \quad
\int_{PD^{\times}_F\backslash PD^{\times}_{\A}} dh = 1. \]
Since $\mathcal{W}_{\psi}(f)$ is nonzero for some $f$, it follows by a standard argument as in [GS, Pg. 2718-2719] that $\mathcal{S}_D(\theta(\phi,f)) \ne 0$ for some $\phi$ and $f$.
The proposition is proved.
\end{proof}
\vskip 10pt

The cuspidal representations $\pi_D$ whose Jacquet-Langlands lifts are not cuspidal are precisely the CAP representations of $GL_2(D)$. Here, the notion of CAP is as given in [G, \S 3.9]. Props.
\ref{P:residue} and \ref{P:shal} essentially show that all CAP representations (as well as the residual spectrum) of $GL_2(D)$ can be obtained as theta lifts from $GL_2$. More precisely, for any such $\pi_D$, there is a unique cuspidal representation $\tau$ of $GL_2$ whose theta lift to $GSO(V_D)$ is
equal to $\pi_D \boxtimes \omega_{\tau}$.

\vskip 15pt
\section{\bf The Local Problem} \label{S:local}
 In this section, we shall study the local analog of Theorems \ref{T:JM2}  and \ref{T:converse}  so as to clarify the nature of the local obstructions there. In particular, we shall relate it to the local Gross-Prasad conjecture. Thus, in this section, we shall let $F$ denote a {\em non-archimedean} local field and $D$ the unique quaternion division algebra over $F$.
\vskip 10pt

Let $\pi_D$ be an irreducible representation of $GL_2(D)$ with
central character $\mu^2$ and let $\pi$ be its Jacquet-Langlands
lift on $GL_4(F)$. Recall that $\pi_D$ (and similarly $\pi$) has
Shalika period with respect to $\mu$ if
\[ \Hom_{\Delta D^{\times} \cdot N_D} (\pi_D,  (\mu \circ \N_D) \boxtimes \psi_D) \ne 0, \]
or equivalently, regarding $\pi_D \boxtimes \mu$ as a representation
of $GSO(V_D)$,
\[  \Hom_{S_D} (\pi_D \boxtimes \mu, \psi_D) \ne 0. \]
It is known that the dimension of this Hom space is at most 1.
\vskip 10pt

In the papers [P1] and [P2], D. Prasad has studied the question of
existence of local Shalika periods, especially for irreducible
principal series representations. Let us recall his results briefly.
Recall that if $\tau_1$ and $\tau_2$ are two infinite-dimensional
representations of $GL_2(F)$, then $PS(\tau_1,\tau_2)$ denotes the
representation of $GL_4(F)$ unitarily induced from the
representation $\tau_1 \boxtimes \tau_2$ of $P$.  Similarly, if
$\tau_{D,1}$ and $\tau_{D,2}$ are two representations of
$D^{\times}$, then $PS(\tau_{D,1}, \tau_{D,2})$ is the analogous
principal series representation induced from $P_D$. The following
proposition is due to D. Prasad [P1, Prop. 7] and [P2, Thm. 2]:
\vskip 10pt

\begin{Prop} \label{P:prasad}
(i) We have a short exact sequence of $GL_2$-modules:
\[  \begin{CD}
0 @>>> \tau_1 \otimes \tau_2 @>>> PS(\tau_1,\tau_2)_{N,\psi} @>>>
\pi(\omega_1 |-|^{1/2}, \omega_2 |-|^{-1/2}) @>>> 0, \end{CD} \]
where $\omega_i$ is the central character of $\tau_i$. Thus, if
$PS(\tau_1,\tau_2)$ is irreducible, then it possesses local Shalika
period with respect to $\mu$ if and only if one of the following
holds:
\begin{itemize}
\item $\omega_1 = \omega_2 = \mu$;
\vskip 5pt
\item $\tau_1^{\vee} = \tau_2 \otimes \mu^{-1}$.
\end{itemize}
\vskip 10pt

(ii) We have an isomorphism of $D^{\times}$-modules:
\[  PS(\tau_{D,1}, \tau_{D,2})_{N_D,\psi_D} \cong \tau_{D,1} \otimes \tau_{D,2}. \]
Thus, if $PS(\tau_{D,1}, \tau_{D,2})$ is irreducible, it possesses
local Shalika period with respect to $\mu$ if and only if
$\tau_{D,1}^{\vee} \cong \tau_{D,2} \otimes \mu^{-1}$.
\end{Prop}
\vskip 10pt

\begin{proof}
The reader may notice that the statement of (i) is different from
that in [P2, Thm. 2]; namely, in the 3rd term of the short exact
sequence,  [P2, Thm. 2] has $\pi(\omega_1,\omega_2)$ instead of
$\pi(\omega_1 |-|^{1/2}, \omega_2 |-|^{-1/2})$. Prasad has informed
us that there is a normalization error in [P2, Thm. 2] and since the
proof given there is somewhat sketchy, he has kindly provided us
with a detailed proof, which we reproduce here. \vskip 5pt

(i) We have:
\[  GL_4 = P \cup Pw_{23}P \cup PwP \]
where $w_{23} = (23) \in S_4$ (the Weyl group of $GL_4$) and $w =
(13)(24)$. By Mackey theory, the restriction of $PS(\tau_1,\tau_2)$
to $P$ has a filtration with successive quotients
\[  \begin{cases}
A = \delta_P^{1/2} \cdot (\tau_1 \boxtimes \tau_2) \\
B =  ind^P_{P \cap w_{23}Pw_{23}} \delta_P^{1/2} \cdot (\tau_1 \boxtimes \tau_2) \\
C = ind^P_{P \cap wPw} \delta_P^{1/2} \cdot (\tau_1 \boxtimes
\tau_2).  \end{cases} \] We shall be interested in the restriction
of these representations to the Shalika subgroup $\tilde{S}$. Now
the quotient $A$ does not contribute to the twisted Jacquet module
since $N$ acts trivially. As for $C$,  $P \cap wPw = M = GL_2 \times
GL_2$ and one sees that
\[  C_{N, \psi} \cong \tau_1 \otimes \tau_2 \]
as representations of $\Delta GL_2$. Thus it remains to show that
\[  B_{N,\psi} \cong  \pi(\omega_1 |-|^{1/2}, \omega_2 |-|^{-1/2}). \]
\vskip 5pt

The stabilizer group $P \cap w_{23}Pw_{23}$ is equal to $(B \times
B) \cdot N_0$, where
\[  N_0 = \{ \left( \begin{array}{cc}
* & * \\
0 & * \end{array} \right)  \} \subset N. \] Now consider the
restriction of $B$ to the Shalika subgroup $\tilde{S} = \Delta GL_2
\cdot N$. The double coset space
\[  \tilde{S} \backslash P / (P \cap PwP) = \Delta GL_2 \backslash (GL_2 \times GL_2)/(B \times B) \]
has size $2$. By Macket theory, we  need to consider each of these
double cosets in turn. \vskip 5pt

Let us first consider the non-trivial double coset which is
represented by the Weyl group element $w_{12} = (12)$. We first
compute:
\[ \Delta T \cdot N^0 \]
where
\[  N^0 =  \{ \left( \begin{array}{cc}
0 & * \\
* & * \end{array} \right)  \} \subset N. \]
So the representation under consideration is
\[  Ind^{\Delta GL_2 \cdot N}_{\Delta T \cdot N^0} \tau_1 \otimes \tau_2, \]
where the action of an element in $T \cdot N^0$ is via the sequence
of maps
\[  t \cdot n \mapsto w_{12} (t \cdot n) w_{12} \mapsto w_{23}w_{12}(t \cdot n) w_{12} w_{23} \in P \]
followed by the action of $P$ on $\tau_1 \boxtimes \tau_2$.
Explicitly,
\[
\left(
\begin{array}{cccc}
t_1 &  & 0 & n_3 \\
 & t_2 & n_1 & n_2 \\
 &  &  t_1 &  \\
 & &  &  t_2  \end{array} \right)   \mapsto
\left(  \begin{array}{cccc}
t_2 & n_1 & 0 & n_2 \\
 & t_1 & 0 & 0\\
 &  &  t_1 & n_3 \\
 & &  &  t_2  \end{array} \right). \]
Since the character $\psi$ of $N$ is non-trivial on the 1-parameter
subgroup with coordinate $n_2$, but the restriction of $\tau_1
\boxtimes \tau_2$ to the image of this 1-parameter subgroup under
the above map is trivial, we see that the non-trivial double coset
does not contribute to the twisted Jacquet module. \vskip 5pt

It remains to consider the trivial double coset which gives rise to
the representation
\[  Ind^{\Delta GL_2 \cdot N}_{\Delta B \cdot N_0} \delta_P^{1/2} \cdot (\tau_1 \otimes \tau_2). \]
We note the following general lemma: \vskip 5pt

\begin{Lem}
For a representation $\pi$ of $\Delta B \cdot N$,
\[ \left(  Ind^{\Delta GL_2 \cdot N}_{\Delta B \cdot N} \pi \right)_{N,\psi} \cong Ind^{\Delta GL_2}_{\Delta B}
\pi_{N, \psi} \]
\end{Lem}
\vskip 5pt

Therefore, we need to calculate
\[  \left( ind_{\Delta B \cdot N_0}^{\Delta B \cdot N} \tau_1 \otimes \tau_2 \right)_{N,\psi} \]
as a $\Delta B$-module. For this, note the following: \vskip 5pt

\begin{Lem}
Let $\psi_0 = \psi|_{N_0}$. Then for a representation $\pi$ of
$\Delta B \cdot N_0$,
\[  \left( ind_{\Delta B \cdot N_0}^{\Delta B \cdot N} \pi \right)_{N,\psi} = \delta_B^{-1} \cdot \pi_{N_0,\psi_0}  \]
where $\delta_B$ is the modulus character of $\Delta B$ (which is
the inverse of the character of the action of $\Delta B$ on
$N/N_0$).
\end{Lem}

\vskip 5pt

Applying this lemma, and using the fact that the action of $\Delta B
\cdot N_0$ on $\delta_P^{1/2} \cdot (\tau_1 \otimes \tau_2$) is via
the map
 \[  \left(
\begin{array}{cccc}
a & b & n_1 & n_2 \\
 & c & 0 & n_3 \\
 &  &  a & b \\
 & &  &  c  \end{array} \right)   \mapsto
\left(  \begin{array}{cccc}
a & n_1 & b & n_2 \\
 & a & 0 & b\\
 &  &  c & n_3 \\
 & &  &  c  \end{array} \right), \]
we conclude that
\[  (\delta_p^{1/2} \cdot (\tau_1 \otimes\tau_2))_{N_0,\psi_0} = \delta_B \cdot (\omega_1 \boxtimes \omega_2) \]
as representations of $\Delta B$. Putting everything together
completes the proof of (i). \vskip 5pt

(ii) The proof is similar and in fact easier; we refer the reader to
[P1, Prop. 7].
\end{proof}

\vskip 10pt From this proposition, one sees that it is possible that
a representation $\pi$ of $GL_4$ has local Shalika period, but its
Jacquet-Langlands lift $\pi_D$ does not. This was exploited in our
construction of  the counterexample to Thm. \ref{T:JM2}. Note
however that our local condition in Thm. \ref{T:converse} rules out
more representations $\pi$ than those for which $\pi_D$ has no
Shalika period. For example, if $\tau_1 = \tau_2$ has central
character $\mu$, then by the Proposition, both $PS(\tau_1,\tau_2)$
and its Jacquet-Laglands transfer to $GL_2(D)$ admit Shalika
period with respect to $\mu$.  However, such representations are
ruled out by the local condition in Thm. \ref{T:converse}. This is
explained by Thm. \ref{T:non-split}. \vskip 10pt

To study the local Shalika period of more general representations,
such as the discrete series, we note the following local analog of
Proposition \ref{P:main}. \vskip 5pt

\begin{Prop}
(i) As representations of $GSO(V_D)$,
\[  (\Omega_D)_{U, \chi} = ind_{S_D}^{GSO(V_D)} \psi_D. \]
\vskip 5pt

(ii) Suppose that  $D$ is split and $U_D$ is a maximal unipotent
subgroup of $GSO(V_D)$ with nondegenerate character $\chi_D$. Then
as representations of $GSp_4$,
\[ (\Omega_D)_{U_D, \chi_D} \cong ind_U^{GSp_4} \chi. \]
\end{Prop}

\begin{proof}
The statement (i) is proved in a similar way as Prop. \ref{P:main}.
The statement (ii) is essentially [MS, Prop. 4.1]; see also [GRS,
Prop. 2.4 and Cor. 2.5].
\end{proof}
\vskip 10pt

\begin{Cor}  \label{C:local}
(i) Let $\pi_D$ be an irreducible  representation of $GSO(V_D)$ with
central character $\mu^2$. Then  $\pi_D$ has Shalika period with
respect to $\mu$  if and only if $\dim \Theta(\pi_D \boxtimes
\mu)_{U,\chi}=1$. \vskip 5pt

(ii) Suppose that  $D$ is split and $\pi_D$ is generic. Then
$\pi_D$ has Shalika period with respect to $\mu$ if and only if
$\theta(\pi_D \boxtimes \mu)$ is generic.  Moreover, if $\sigma$
is an irreducible representation of $GSp_4$, then $\sigma$ is
generic if and only if $\Theta_D(\sigma)$ is generic, in which case
$\dim \Theta_D(\sigma)_{U_D,\chi_D} =1$.
\end{Cor}

\vskip 10pt

\begin{Thm}  \label{T:split}
Let $\pi$ be a generic representation of $GL_4$ with central
character $\mu^2$. The following are equivalent: \vskip 5pt

\noindent (i) $\pi$ has Shalika period with respect to $\mu$; \vskip
5pt

\noindent (ii) the small theta lift $\theta(\pi\boxtimes \mu)$ of
$\pi\boxtimes \mu$ to $GSp_4$ is generic; \vskip 5pt

\noindent (iii)  the small theta lift $\theta(\pi\boxtimes \mu)$ of
$\pi\boxtimes \mu$ to $GSp_4$ is non-zero; \vskip 5pt

\noindent (iv) the Langlands parameter $\varphi_{\pi}$ of $\pi$
factors through $GSp_4(\C)$ and its similitude character $\Lambda \circ \varphi_{\pi}$ is equal to $\mu$.
\vskip 5pt

If $\pi$ is a discrete series representation, then the above
conditions are equivalent to: \vskip 5pt

\noindent (v) The $L$-factor $L(s, \pi, \bigwedge^2 \otimes \mu^{-1})$ defined by Shahidi has a pole
at $s = 0$.
\end{Thm}

\begin{proof}
The equivalence of (i), (ii) and (iii) is the previous corollary.
 For the equivalence of (iv) and (v) in the case of discrete series representations, note that if
\[ \varphi_{\pi}: W_F' \longrightarrow GL_4(\C) \]
is the Langlands parameter of $\pi$,  then a recent result of Henniart [He] shows that the local Langlands correspondence for $GL_n$ respects twisted exterior square L-functions, so that
\[  L(s,\pi, \bigwedge^2 \otimes \mu^{-1}) = L(s, \bigwedge^2 \varphi_{\pi} \otimes \mu^{-1}). \]
Moreover,  $L(s, \bigwedge^2 \varphi_{\pi} \otimes \mu^{-1})$  has a pole at $s = 0$ if and only if $\left(
\bigwedge^2 \varphi_{\pi} \right) \otimes \mu^{-1}$ contains the
trivial representation as a summand. In other words, the action of
$W_F'$ via $\varphi_{\pi}$ preserves a non-zero symplectic form up
to scaling by the character $\mu$ (thought of as a character of
$W'_F$ by local class field theory). This symplectic form is necessarily
nondegenerate, so that $\varphi_{\pi}$ factors through $GSp_4(\C)$, for
otherwise, its kernel is a non-trivial $W'_F$-submodule, which contradicts the irreducibility of
$\varphi_{\pi}$.
  \vskip 10pt

 The main assertion of the theorem is thus the equivalence of (i)-(iii) and (iv). In fact, the key case of discrete series representations is a special  case of a beautiful theorem of Muic-Savin [MS, Thm. 2.2], which shows that the theta lift of a discrete series representation $\pi$ to $GSp_4$ is non-zero iff
 (v) holds.
 \vskip 5pt

 With the discrete series case taken care of, a non-discrete series generic representation $\pi$ is of the form $Ind_Q^{GL_4} \tau$ with $\tau$ a twist of a discrete series representation on the Levi factor of some parabolic $Q$. If $Q$ is not the $(1,3)-$ or $(3,1)-$ parabolic, then by induction-in-stages, we may assume that $Q$ is the parabolic $P$, in which case $\tau$ is generic but not necessarily a discrete series. The equivalence of (i) and (iv) then follows readily from Prop. \ref{P:prasad}.
 \vskip 5pt

 It remains to consider the case when $\pi = Ind_Q^{GL_4} \tau$ with $Q$ the $(3,1)$-parabolic and
 $\tau$ is the a twist of a discrete series representation. In this case, the Langlands parameter of $\pi$ is not of symplectic type and so we need to show that $\pi$ does not have local Shalika period with respect to $\mu$. This can be checked by a Mackey theory argument analogous to the proof of Prop.
 \ref{P:prasad}.
\end{proof}

\vskip 5pt

\noindent{\bf Remarks:} In the recent preprint of Jiang-Nien-Qin [JNQ], the case of a general division algebra $D$ of degree $\geq 2$ is considered and the equivalence of (i), (iv) and (v) for supercuspidal
$\pi$ and trivial $\mu$ is shown.
\vskip 10pt

\vskip 10pt Since it is not much of a trouble and for the
convenience of reference, we list the different types of symplectic parameters
$\varphi_{\pi} : W_F' \longrightarrow GSp_4(\C)$ with similitude
character $\mu$ which could give rise to generic representations of $GL_4$:
\vskip 5pt

\noindent (1) $\varphi_{\pi}$ is irreducible. In this case, $\pi$ is
a discrete series representation. \vskip 5pt

\noindent (2) $\varphi_{\pi} = \phi_1 \oplus \phi_2$ where each
$\phi_i$ is an irreducible 2-dim representation satisfying one of
the following: \vskip 5pt (a) $\det \phi_i = \mu$; \vskip 5pt (b)
$\phi_1^{\vee} = \phi_2 \otimes \mu^{-1}$. \vskip 5pt In this case,
$\pi$ is the representation $PS(\tau_1,\tau_2)$ where $\tau_i$ is
the discrete series representation of $GL_2$ associated to $\phi_i$.
\vskip 5pt

\noindent (3) $\varphi_{\pi} = \phi \oplus \chi_1 \oplus \chi_2$,
where $\phi$ is irreducible of dimension 2 and $\det \phi = \mu$,
whereas the $\chi_i$'s are 1-dimensional with $\chi_1 \chi_2  =\mu$.
In this case, $\pi = PS(\tau, \pi(\chi_1,\chi_2))$ where $\tau$ is
the discrete series representation associated to $\phi$. \vskip 5pt

\noindent (4) $\varphi_{\pi} = \chi_1 \oplus \mu \chi_1^{-1} \oplus
\chi_2 \oplus \mu\chi_2^{-1}$. In this case, $\pi = PS(\pi(\chi_1,
\mu \chi_1^{-1}), \pi(\chi_2, \mu \chi_2^{-1}))$. \vskip 10pt

Of these 4 classes, only (1) and (2) are relevant parameters for
$GL_2(D)$. Moreover, note that the cases (2a) and (2b) are not
disjoint: their intersection consists of those parameters with
$\phi_1 = \phi_2$ and $\det \phi_i = \mu$.

\vskip 10pt

Now we examine the analogous problem for $\pi_D$.  This may be a
good time to bring in the local Gross-Prasad conjecture [GP].  We
have been considering Shalika priods with respect to $1$-dimensional
representations of the Shalika group, by using the character $\mu$
on $D^{\times}$ or $GL_2$. More generally,   one may consider the
generalized Shalika period with respect to  any irreducible
representation $\tau$ of $D^{\times}$ or $GL_2$. For example, on
$GL_4$, one may consider the generalized Shalika period with respect
to the twisted Steinberg representation $St \otimes \mu$. The local
Gross-Prasad conjecture predicts that given tempered representations
$\pi_D$ and $\pi = JL(\pi_D)$,
\[  \dim \Hom_{\tilde{S}_D}(\pi_D, \mu \boxtimes \psi_D) + \dim \Hom_{\tilde{S}}(\pi, (St\otimes \mu) \boxtimes \psi) =1, \]
and
\[   \dim \Hom_{\tilde{S}_D}(\pi_D, \mu \boxtimes \psi_D) = 1 \Longleftrightarrow
\epsilon( (\bigwedge^2 \varphi_{\pi} \otimes \mu^{-1}) \otimes S_2)= -1, \]
where $S_2$ denotes the 2-dimensional representation of $W_F'$ trivial on $W_F$ (so essentially it is a representation of $SL_2(\C)$) and is the Langlands parameter of the Steinberg representation $St$.
This last equivalence is part of Theorem \ref{T:non-split} below. At
this point, we note the following lemma which is explained to us by
D. Prasad: \vskip 5pt

\begin{Lem} \label{L:epsilon}
Suppose that $\pi$ has central character $\mu^2$. Then the following
are equivalent:
\vskip 5pt

(i)  $\epsilon( (\bigwedge^2 \varphi_{\pi} \otimes \mu^{-1}) \otimes S_2) = -1$.
\vskip 5pt

(ii) $\dim((\bigwedge^2 \varphi_{\pi} \otimes \mu^{-1})^{W'_F})$ is
odd. \vskip 5pt

(iii) $\varphi_{\pi}$ is of type (1) or (2b). \vskip 5pt

\noindent In particular, if $\epsilon( (\bigwedge^2 \varphi_{\pi} \otimes \mu^{-1}) \otimes S_2) = -1$, then $L(s, \bigwedge^2 \varphi_{\pi} \otimes \mu^{-1})$ has a pole at $s= 0$.
\end{Lem}
\vskip 5pt

\begin{proof}
The representation $\bigwedge^2 \varphi_{\pi} \otimes \mu^{-1}$ is a
map
\[  W'_F \longrightarrow SO_6(\C). \]
In particular, it is a self dual representation with trivial
determinant. We shall calculate the $\epsilon$-factor by considering
different cases. The reader can find similar computations in [P2].
The following well-known identity, for which we refer to [P2, \S 5],
will be very useful:
\[
\epsilon(\rho \otimes S_n,\psi) = \epsilon(\rho,\psi)^n \cdot  \det
(-Frob| \rho^I)^{n-1},
\]
where $\rho$ is a representation of the Weil-Deligne group $W_F'$ trivial on $SL_2(\C)$ and $S_n$
is the $n-$dimensional irreducible representation of $SL_2(\C)$
regarded as a representation of $W'_F$. In
addition, we shall use the fact that
\[  \epsilon(\rho,\psi) \cdot \epsilon(\rho^*,\psi) = \det \rho(-1). \]
In particular, if $\rho^* = \rho$, then $\epsilon(\rho,\psi)^4 =1$
and $\epsilon(\rho,\psi)^2=1$ if $\det \rho$ is trivial.
Now we may consider the different cases, according to how $\varphi_{\pi}$ decomposes as a representation of $SL_2(\C)$.

\vskip 15pt

\noindent \underline{{\bf Case 1:} $\varphi_{\pi}$ is a
representation of the Weil group $W_F$} \vskip 5pt

Applying the above identities for the epsilon factor, we see that
\[  \epsilon(S_2 \otimes (\bigwedge^2 \varphi_{\pi} \otimes \mu^{-1}))  =
  \det (-Frob| (\bigwedge^2 \varphi_{\pi} \otimes \mu^{-1})^I). \]
Here, if  $(\bigwedge^2 \varphi_{\pi} \otimes \mu^{-1})^I$ is zero,
then the determinant in question is interpreted to be $1$. Now the
space $(\bigwedge^2 \varphi_{\pi} \otimes \mu^{-1})^I$ is precisely
the submodule spanned by the unramified characters occurring in
$\bigwedge^2 \varphi_{\pi} \otimes \mu^{-1}$. By self-duality, if an
unramified character $\chi$ occurs, then so must its inverse
$\chi^{-1}$. If $\chi^{-1} \ne \chi$, then the determinant of
$-Frob$ on this 2-dimensional submodule is 1. On the other hand, if
$\chi^{-1} = \chi$, then $\chi$ is either the trivial character or
the unique unramified quadratic character. The action of $-Frob$ on
$\chi$ is then $-1$ and $1$ respectively. So we see that:
\[   \det (-Frob| (\bigwedge^2 \varphi_{\pi} \otimes \mu^{-1})^I) =
(-1)^{\dim(\bigwedge^2 \varphi_{\pi} \otimes \mu^{-1})^{W'_F}}. \]
 This shows the equivalence of (i) and (ii) in this case. Moreover, if $\varphi_{\pi}$ is irreducible, then
 \[  \epsilon(S_2 \otimes (\bigwedge^2 \varphi_{\pi} \otimes \mu^{-1}))  =-1 \]
 iff  $\varphi_{\pi}$ is of symplectic type with similitude character $\mu$.
 \vskip 10pt

\noindent \underline{{\bf Case 2:} $\varphi_{\pi} = \rho \oplus
(\chi \cdot S_2)$} \vskip 5pt

In this case, $\det(\rho) \cdot \chi^2 = \mu^2$ and
\[  \bigwedge^2\varphi_{\pi} \otimes \mu^{-1} = \mu^{-1} \cdot \det (\rho) \oplus \mu^{-1} \cdot \chi^2 \oplus \mu^{-1}\cdot  \chi \cdot (\rho \otimes S_2) \]
and so
\[
S_2 \otimes  \bigwedge^2\varphi_{\pi} \otimes \mu^{-1} =
((\mu^{-1}\det(\rho) \oplus \mu^{-1} \cdot \chi^2) \otimes S_2)
\oplus \mu^{-1} \chi \cdot \rho \oplus \mu^{-1}\chi \cdot (\rho
\otimes S_3). \] Observe that $(\mu^{-1} \cdot \det(\rho) \oplus
\mu^{-1} \cdot \chi^2)$ and $\mu^{-1} \chi \cdot \rho$ are both
self-dual with determinant $1$. Thus, we see that
\[  \epsilon(S_2 \otimes (\bigwedge^2 \varphi_{\pi} \otimes \mu^{-1})) =  \det(-Frob| (\mu^{-1}\det(\rho))^I) \cdot \det(-Frob|(\mu^{-1}\chi^2)^I) = (-1)^{\dim(\bigwedge^2 \varphi_{\pi} \otimes \mu^{-1})^{W'_F}}, \]
which shows the equivalence of (i) and (ii) in this case. Note that
for each of these determinants to be $-1$, we need $\det(\rho) =
\mu$ and $\chi^2 = \mu$ respectively. But since $\det(\rho) \cdot
\chi^2 = \mu^2$, if one of these holds, so does the other. Thus, we
see that the $\epsilon$-factor is always $1$. \vskip 5pt

\noindent \underline{{\bf Case 3:} $\varphi_{\pi} = \rho \otimes
S_2$}

\vskip 5pt

In this case, $\det(\rho)^2 = \mu^2$, so that $\det(\rho) \cdot
\mu^{-1}$ is a quadratic character. We have:
\[  \bigwedge^2\varphi_{\pi} \otimes \mu^{-1}
= \mu^{-1} \cdot  Sym^2 \rho \oplus \mu^{-1} \cdot \det(\rho)
\otimes S_3 \] and
\[  S_2 \otimes \bigwedge^2\varphi_{\pi} \otimes \mu^{-1}
= (\mu^{-1} \cdot Sym^2\rho \otimes S_2)  \oplus
(\mu^{-1}\det(\rho)) \otimes S_2 ) \oplus \mu^{-1}\det(\rho) \otimes
S_4. \] Now observe that $\mu^{-1} \cdot Sym^2(\rho)$ and $\mu^{-1}
\det(\rho)$ are self-dual with determinant $\mu^{-1} \det(\rho)$.
Thus, a short calculation gives:
\[  \epsilon(S_2 \otimes (\bigwedge^2 \varphi_{\pi} \otimes \mu^{-1})) =  \det(-Frob| (\mu^{-1} \cdot Sym^2\rho)^I) =  (-1)^{\dim(\bigwedge^2 \varphi_{\pi} \otimes \mu^{-1})^{W'_F}}. \]
This shows the equivalence of (i) and (ii) in this case. Moreover,
if $\rho$ is irreducible, then $\epsilon(S_2 \otimes (\bigwedge^2
\varphi_{\pi} \otimes \mu^{-1})) = -1$ iff $\rho$ is induced from a
character of $W_K$ with $K/F$ a quadratic extension of $F$ and
$\mu^{-1} \cdot \det \rho = \omega_{K/F}$ (the quadratic character
associated to $K/F$ by local class field theory). \vskip 5pt

\noindent \underline{{\bf Case 4:} $\varphi_{\pi} = \chi_1 \oplus
(\chi_2 \otimes S_3)$} \vskip 5pt

We have $\chi_1 \cdot \chi_2^3 = \mu^2$. Now
\[  \bigwedge^2\varphi_{\pi} \otimes \mu^{-1} = (\chi_1\chi_2\mu^{-1} \otimes S_3) \oplus (\chi_2^2 \mu^{-1} \otimes S_3), \]
which contains no trivial representation, and a short calculation
shows that $\epsilon(S_2 \otimes (\bigwedge^2 \varphi_{\pi} \otimes
\mu^{-1}))$ is always equal to $1$.

\vskip 5pt

\noindent \underline{{\bf Case 5:} $\varphi_{\pi} = \chi \otimes
S_4$}

\vskip 5pt We have $\chi^4 = \mu^2$ and
\[  \bigwedge^2\varphi_{\pi} \otimes \mu^{-1} = \chi^2\mu^{-1} \oplus (\chi^2\mu^{-1} \otimes S_5) \]
and a short computation gives
\[ \epsilon(S_2 \otimes (\bigwedge^2 \varphi_{\pi} \otimes \mu^{-1})) = \det(-Frob| (\chi^2\mu^{-1})^I)
   = (-1)^{\dim(\bigwedge^2 \varphi_{\pi} \otimes \mu^{-1})^{W'_F}}, \]
   which shows the equivalence of (i) and (ii). Indeed, the $\epsilon$-factor is $-1$ iff $\chi^2 = \mu$.
   \vskip 10pt

 We have thus shown the equivalence of (i) and (ii) in general. From this, it follows that if the epsilon factor is $-1$, then the trivial representation occurs in $\bigwedge^2 \varphi_{\pi} \otimes \mu^{-1}$
so that $\varphi_{\pi}$ is of symplectic type and $L(s, \pi, \bigwedge^2 \otimes \mu^{-1})$ has a pole at $s= 0$.   A short computation now gives the
following table:

\vskip 10pt
\begin{tabular}{|c|c|c|c|c|}
\hline
& & & & \\
Type of $\varphi_{\pi}$ & (1) & (2b) and (2a) & (2b) but not (2a) & (2a) but not (2b)  \\
& & & & \\
\hline
& & & & \\
$\dim(\bigwedge^2 \varphi_{\pi} \otimes \mu^{-1})^{W'_F}$ & 1 & 3 & 1 & 2 \\
& & & & \\
\hline
\end{tabular}
 \vskip 15pt
From this, we see the equivalence of (ii) and (iii). The lemma is proved.
\end{proof}

Our main local theorem is: \vskip 5pt

\begin{Thm} \label{T:non-split}
Suppose that $\pi_D$ is a representation of $GL_2(D)$ with central
character $\mu^2$ such that its Jacquet-Langlands lift $\pi$ to
$GL_4$ is generic. The following are equivalent: \vskip 5pt

\noindent (i) $\pi_D$ has Shalika period with respect to $\mu$;
\vskip 5pt

\noindent (ii) the big theta lift $\Theta(\pi_D\boxtimes \mu)$ of
$\pi_D\boxtimes \mu$ to $GSp_4$ is generic (and thus non-zero).
\vskip 5pt

\noindent (iii) $\epsilon((\bigwedge^2 \varphi_{\pi} \otimes
\mu^{-1}) \otimes S_2) = -1$.

\vskip 10pt

Moreover, when these conditions hold, the small theta lift
$\theta(\pi_D \boxtimes \mu)$ is non-generic precisely when $\pi_D
= PS(\tau_{D,1}, \tau_{D,2})$ where $\tau_{D,1} = \tau_{D,2}$  has
central character $\mu$, i.e. when $\varphi_{\pi}$ is of both type
(2a) and (2b), or equivalently when $\dim(\bigwedge^2 \varphi_{\pi}
\otimes \mu^{-1})^{W'_F} >  1$.
\end{Thm}

\vskip 5pt

The rest of the section is devoted to the proof of this theorem. The
equivalence of (i) and (ii) is Cor. \ref{C:local}. The content of
the theorem is thus the equivalence of (i)+(ii) and (iii). If $\pi_D
= PS(\tau_{D,1}, \tau_{D,2})$, then the equivalence of these follows
immediately from Prop. \ref{P:prasad}. Moreover, from the explicit
determination of local theta correspondence given in the next
section (specifically Thm. \ref{T:local-theta-D}(i)),
one sees that if $\pi_D = PS(\tau_{D,1},\tau_{D,2})$ has
local Shalika period with respect to $\mu$, then $\theta(\pi_D
\boxtimes \mu)$ is non-generic iff $\tau_{D,1}=\tau_{D,2}$ has
central character $\mu$.

\vskip 5pt

Thus it remains to consider the case of discrete series
representations. We note first that if $\pi_D$ is supercuspidal and
has local Shalika period with respect to $\mu$, then $\theta(\pi_D
\boxtimes \mu) = \Theta(\pi_D \boxtimes \mu)$ is generic. In fact,
the same result holds if $\pi_D$ is a discrete series
representation, by general results of Muic [Mu, Thm. 6.2]. Moreover,
by the previous lemma, it is clear that $\epsilon(
(\bigwedge^2 \varphi_{\pi} \otimes \mu^{-1}) \otimes S_2) = -1$ if and only if
$\varphi_{\pi}$ is of symplectic type with similitude character
$\mu$. Hence, to complete the proof of Thm. \ref{T:non-split}, it
remains to show: \vskip 10pt

\begin{Thm} \label{T:local-JM}
Let $\pi_D$ be a discrete series representation of $GL_2(D)$ with
central character $\mu^2$ and Jacquet-Langlands lift $\pi$ on
$GL_4$. Then $\pi_D$ has Shalika period with respect to $\mu$ if and
only if $\pi$ does.
\end{Thm}
\vskip 5pt

It is clear that the theorem follows from the following proposition,
which addresses the very natural question about the compatibility of
the theta correspondence and the Jacquet-Langlands transfer of
discrete series representations. \vskip 5pt

\begin{Prop} \label{P:JL}
Suppose that $\pi_D$ is a discrete series representation with
Jacquet-Langlands lift $\pi$ on $GL_4$. \vskip 5pt

\noindent (i) If $\pi_D$ has local Shalika period with respect to
$\mu$, so that $\sigma_D  = \theta(\pi_D \boxtimes \mu)$ is
generic, then the small theta lift of $\sigma_D$ to $GSO(V)$ is
isomorphic to $\pi \boxtimes \mu$. \vskip 5pt

\noindent (ii) Conversely, if $\pi$ has local Shalika period with
respect to $\mu$, so that $\sigma = \theta(\pi \boxtimes \mu)$ is
generic, then the small theta lift of $\sigma$ to $GSO(V_D)$ is
isomorphic to $\pi_D$.
\end{Prop}

\begin{proof}
Our proof of the proposition is going to involve global arguments.
Let us consider the statement (ii), so that we are starting with a
discrete series $\pi \boxtimes \mu$ on $GSO(V)$ with Shalika period.
By Theorem \ref{T:split}, we see that in the terminology of the
proof of Lemma \ref{L:epsilon}, $\pi$ is one of the following:
\vskip 5pt

\begin{itemize}
\item $\pi$ is supercuspidal with symplectic parameter and similitude character $\mu$;
\vskip 5pt
\item $\pi$ is a generalized Steinberg representation attached to a dihedral supercuspidal representation $\pi_{\rho}$, i.e. $\pi$ has parameter of the form $\rho \boxtimes S_2$ (as in Case 3 in the proof of Lemma \ref{L:epsilon}) with $\rho$ irreducible and monomial with respect to a quadratic extension $K/F$ and $\det \rho = \mu \cdot \omega_{K/F}$;
\vskip 5pt
\item $\pi$ is a twisted Steinberg $St_{\chi}$ with $\chi^2 = \mu$ (as in Case 5 in the proof of Lemma \ref{L:epsilon}).
\end{itemize}
\vskip 5pt

The main technical tool we need is: \vskip 5pt

\begin{Lem}
Let $\pi \boxtimes \mu$ be as in (ii) of the proposition.  Let
$\mathbb{F}$ be a number field such that $\mathbb{F}_v = F$ for some
place $v$ of $\mathbb{F}$. Then
 $\pi$ can be globalized to a cuspidal representation $\Pi$ of $GL_4$ over $\mathbb{F}$ such that
$\Pi$ has global Shalika period with respect to some $\Upsilon$ such
that $\Upsilon_v = \mu$.
\end{Lem}

\begin{proof}
If $\pi$ is supercuspidal, this is a consequence of a general result
of Prasad and Schulze-Pillot [PSP, Thm. 3.1] (proved using a simple
form of the relative trace formula!). So suppose that $\pi$ is not
supercuspidal, so that it is of the other two types described above.
If $\pi = St_{\chi}$ with $\chi^2 = \mu$, then note that $St_{\chi}$
is the Langlands lift of the twisted Steinberg representation
$st_{\chi}$ of $GL_2(F)$ under the adjoint cube lifting $Sym^3
\otimes \det^{-1}$. Let $\Omega$ be a cuspidal representation of
$GL_2$ over $\mathbb{F}$ such that $\Omega_v = st_{\chi}$ and the
central character of $\Omega$ is $\Upsilon$. By the results of
Kim-Shahidi [KS], we may consider the adjoint cube lifting of
$\Omega$ to get a cuspidal representation $\Pi$ of $GL_4$ so that
the central character of $\Pi$ is $\Upsilon^2$ and $\Omega_v =
St_{\chi}$. Moreover, by [KS, Pg. 877], the partial twisted exterior
square L-function $L^S(s, \Pi, \bigwedge^2 \otimes \Upsilon^{-1})$
has a pole at $s=1$. This $\Pi$ is the desired cuspidal
representation. \vskip 5pt

Finally, suppose $\pi$ has parameter of the type $\rho \boxtimes
S_2$.  Then observe that $\pi$ is the Langlands lift of the
representation $\pi_{\rho} \boxtimes st$ of $GL_2 \times GL_2$ under
the Rankin-Selberg lifting $GL_2 \boxtimes GL_2 \longrightarrow
GL_4$. We may find a dihedral cuspidal representation $\Omega$ of
$GL_2$ associated to a quadratic extension $\mathbb{K}/\mathbb{F}$
such that $\Omega_v = \pi_{\rho}$ and the central character of
$\Omega$ is $\Upsilon \cdot \omega_{\mathbb{K}/\mathbb{F}}$.
Similarly, let $\Omega'$ be a cuspidal representation of $PGL_2$
such that $\Omega'_v = st$. Then by results of Ramakrishnan [R], the
cuspidal representation $\Omega \boxtimes \Omega'$ has a functorial
lifting to a cuspidal $\Pi$ on $GL_4$ whose local component at $v$
is $\pi$ and such that $L^S(s, \Pi, \bigwedge^2 \otimes
\Upsilon^{-1})$ has a pole at $s=1$. This $\Pi$ is what we are
looking for.
\end{proof}

   \vskip 5pt
As is evident from its proof,  the lemma applies to any finite set
of finite places, though we have stated it only for a singleton set.
To apply the lemma,  let $\mathbb{F}$ be a number field such that
for two places $v_1$ and $v_2$, we have $\mathbb{F}_{v_i} \cong F$
for $i=1$ and $2$.
 Let $\mathbb{D}$ be a global quaternion algebra over $\mathbb{F}$ ramified precisely at $v_1$ and $v_2$. By the lemma, one can find a cuspidal $\Pi$ on $GL_4$ and a character $\Upsilon$ such that
 \vskip 5pt

 \begin{itemize}
 \item $\Pi_{v_i} \cong \pi$ for $i =1$ and $2$;
 \item $\Upsilon_{v_i} = \mu$ for $i =1$ and $2$;
 \item $\Pi$  has global Shalika period with respect to $\Upsilon$.
 \end{itemize}
 \vskip 5pt

Then $\Pi \boxtimes \Upsilon$ has non-zero globally generic cuspidal
theta lift $\Sigma$ on $GSp_4$.  Moreover, by Theorem
\ref{T:converse} (or rather the remark following its proof), we
deduce that $\Sigma$ has non-zero cuspidal theta lift
$\Theta_{\mathbb{D}}(\Sigma)$ to $GSO(V_{\mathbb{D}})$. Now since
   $\Theta_{\mathbb{D}}(\Sigma)$ and the Jacquet-Langlands transfer of $\Pi \boxtimes \Upsilon$ are nearly equivalent, it follows by the strong multiplicity one result of Badulescu \cite[Thm. 5.1(c)]{B}  that they are in fact isomorphic. Extracting the component at $v_1$ proves (ii).
 \vskip 10pt

 For (i), we start with $\pi_D$ on $GL_2(D)$ and get $\sigma_D = \theta(\pi_D \otimes \mu)$ on $GSp_4(F)$. Let $\pi' \boxtimes \mu$ be the theta lift of $\sigma$ to $GSO(V)$ so that $\pi'$ has Shalika period with respect to $\mu$. By (ii), we conclude that $\pi'$ and $\pi_D$ are related by the Jacquet-Langlands correspondence. This proves (i).

 \end{proof}

\vskip 10pt

This completes the proof of Thm. \ref{T:non-split}. It seems to us
that the proof of Prop. \ref{P:JL}  for non-supercuspidal
representations is an overkill, in the sense that it makes use of
too much global machinery. There is in fact a purely local proof of
the proposition in the non-supercuspidal case, by the explicit
determination of the local theta correspondence between $GSp_4$ and
$GSO(V_D)$ (for $D$ both split and non-split). We discuss this in
the next section. \vskip 15pt

\section{\bf Explicit Local Theta Correspondence}\label{S:explicit}

In this section, we describe the explicit determination of local
theta correspondences for the (almost) dual pairs
\[  GSp_4 \times GSO(V) \quad \text{and} \quad GSp_4 \times GSO(V_D), \]
as much as is needed for the applications of this paper.
Specifically, we shall be interested in the theta lift of discrete
series representations in both cases, and the theta lift of
irreducible principal series in the non-split case. As we mentioned
in the introduction, Waldspurger [W] has already determined part of
the correspondence in the split case and our results here are a
refinement of his. \vskip 5pt

To state the results, we introduce some notations. Recall from
Section \ref{S:theta} that we have a Witt decomposition $W = X
\oplus Y$. Suppose that $X = F \cdot e_1 \oplus F \cdot e_2$ and $Y =
F\cdot f_1 \oplus F \cdot f_2$ and consider the decomposition $W = Fe_1
\oplus W' \oplus Ff_1$, where $W' = \langle e_2, f_2 \rangle$. Let
$Q(Z) = L(Z) \cdot U(Z)$ be the maximal parabolic stabilizing the
line $Z = F \cdot f_1$, so that
\[  L(Z) = GL(Z) \times GSp(W') \]
and $U(Z)$ is a Heisenberg group:
\[  \begin{CD}
1 @>>> Sym^2 Z @>>> U(Z) @>>> W' \otimes Z @>>> 1. \end{CD} \] A
representation of $L(Z)$ is thus of the form $\chi \boxtimes \tau$
where $\tau$ is a representation of $GSp(W') \cong GL_2$. We let
$I_{Q(Z)}(\chi,\tau)$ be the corresponding parabolically induced
representation (i.e. via normalized induction).  The module structure of this induced representation
is known.  In particular, we note the following lemma (cf. \cite[Prop. 5.1]{W} and  \cite{ST}):
\vskip 5pt

\begin{Lem}
(a) Let $\tau$ be a supercuspidal representation of $GL_2$. The
induced representation $I_{Q(Z)}(\chi, \tau)$ is reducible iff one
of the following holds: \vskip 5pt

(i) $\chi = 1$; \vskip 5pt

(ii) $\chi = \chi_0 |-|^{\pm 1}$ and $\chi_0$ is a non-trivial
quadratic character  such that $\tau \otimes \chi_0 \cong \tau$.
\vskip 5pt

In case (i), the representation $I_{Q(Z)}(1, \tau)$ is the direct
sum of two irreducible representations, exactly one of which is
generic.  In case (ii), assuming without loss of generality that
$\chi = \chi_0  \cdot |-|$, one has a (non-split) short exact
sequence:
\[  \begin{CD}
0 @>>>   St(\chi_0, \tau_0) @>>> I_{Q(Z)}(\chi_0\cdot |-|, \tau_0 \cdot |-|^{-1/2}) @>>> Sp(\chi_0,\tau_0)
@>>> 0 \end{CD} \]
 where $St(\chi_0,\tau_0)$ is a (generic) discrete series representation and the Langlands quotient $Sp(\chi_0,\tau_0)$ is non-generic.
\vskip 10pt

(b) If $\tau$ is the twisted Steinberg representation  of
$GL_2$, then $I_{Q(Z)}(\chi, \tau)$ is reducible iff one of the following holds:
\vskip 5pt

(i) $\chi = 1$; \vskip 5pt

(ii) $\chi = |-|^{\pm 2}$.
\vskip 5pt

In case (i), $I_{Q(Z)}(1, st_{\chi})$ is the sum of two irreducible representations, exactly one of which is generic. In case (ii), $I_{Q(Z)}(|-|^2, st_{\chi} \cdot |-|^{-1})$ has the
twisted Steinberg representation $St_{PGSp_4} \otimes \chi$ as a
unique irreducible submodule. \vskip 10pt

(c) For general $\tau$, there is a standard intertwining operator
\[ I_{Q(Z)}(\chi^{-1},  \tau \otimes \chi) \longrightarrow I_{Q(Z)}(\chi, \tau). \]
\end{Lem}

\vskip 5pt

Now consider the group $GSO(V_D)$ where $D$ is possibly split. We
may identify $GSO(V_D)$ as a quotient of $GL_2(D) \times GL_1$ as in
Section \ref{S:theta}. We have:
\[  V_D = F \cdot (1,0) \oplus D \oplus F \cdot (0,1) \]
and the stabilizer $P(J)$ of $J = F \cdot (1,0)$ is the image of the
parabolic $P_D \times GL_1 \in GL_2(D) \times GL_1$. A
representation of its Levi subgroup is thus of the form $( \tau_1
\boxtimes \tau_2) \boxtimes \chi$ with $\omega_{\tau_1} \cdot
\omega_{\tau_2} = \chi^2$, and we denote the associated induced
representation of $GL_2(D) \times GL_1$ by $PS(\tau_1,\tau_2)
\boxtimes \chi$. Now we note the following lemma (cf. \cite{T}):
 \vskip 5pt

\begin{Lem}
 Let $\tau$ be a representation of $D^{\times}$ (where $D$ is possibly split)
 and let $JL(\tau)$ denote its Jacquet-Langlands lift to $GL_2$.
 \vskip 5pt

 (i) Suppose that  $JL(\tau)$ is supercuspidal. Then one has  a short exact sequence of
representations of $GL_2(D)$:
\[  \begin{CD}
0 @>>> St(\tau) @>>> PS(\tau |-|^{1/2}, \tau |-|^{-1/2}) @>>>
Sp(\tau) @>>> 0
\end{CD} \]
where $St(\tau)$ is a discrete series representation (a generalized
Steinberg representation) and $Sp(\tau)$ is the unique Langlands
quotient (a generalized Speh representation). \vskip 10pt

(ii) Suppose that $JL(\tau)$ is the twisted Steinberg representation $st_{\chi}$. Then the
principal series $PS(\tau|-|, \tau|-|^{-1})$ has a unique
irreducible submodule which is the twisted Steinberg representation
$St_{\chi}:= St_{GL_2(D)} \otimes \chi$.
\end{Lem}
\vskip 15pt

We can now state the two main theorems of this section: \vskip 5pt

\begin{Thm}  \label{T:local-theta-D}
Consider the case when $D$ is non-split. \vskip 5pt

(i) ({\bf Principal series}) The irreducible principal series
representation $PS(\tau_{D,1}, \tau_{D,2}) \boxtimes \mu$ (with
$\mu^2 = \omega_{D,1} \cdot \omega_{D,2}$) participates in the local
theta correspondence with $GSp_4$ iff one of the following holds:
\vskip 5pt

\begin{itemize}
\item[(a)]  $\omega_{D,1} = \omega_{D,2} = \mu$;

\item[(b)] $\tau_{D,1}^{\vee} \cong \tau_{D,2} \otimes \mu^{-1}$ and  $\omega_{D,1} \ne \mu$.
\end{itemize}

\vskip 5pt Moreover, if (a) holds, then $\theta(PS(\tau_{D,1},
\tau_{D,2}) \boxtimes \mu)$ is the non-generic representation of
$GSp_4(F)$ which is the theta lift of the (supercuspidal)
representation $\tau_{D,1} \boxtimes \tau_{D,2}$ of $GSO(D)$. If
$\tau_{D,1} \ne \tau_{D,2}$, then $\theta(PS(\tau_{D,1}, \tau_{D,2})
\boxtimes \mu)$ is supercuspidal. If $\tau_{D,1} = \tau_{D,2} =
\tau_D$, then $\theta(PS(\tau_{D,1}, \tau_{D,2}) \boxtimes \mu)$ is
the unique non-generic summand of the tempered representation
$I_{Q(Z)}(1, JL(\tau_D))$, which is denoted by $\pi_{ng}(JL(\tau_D))$ in \cite{GT}. \vskip 5pt

If (b) holds,  then
\[  \theta(PS(\tau_{D,1}, \tau_{D,2}) \boxtimes \mu)=I_{Q(Z)}(\frac{\mu}{\omega_{D,1}}, JL(\tau_{D,2}) \cdot \frac{\omega_{D,1}}{\mu})
= I_{Q(Z)}(\frac{\omega_{D,1}}{\mu}, JL(\tau_{D,2})) \] which is
irreducible and generic. \vskip 10pt

(ii) ({\bf Generalized Steinberg})  If $\dim \tau_D > 1$ and $\mu =
\omega_{\tau_D} \cdot \chi$ with $\chi^2 = 1$, then
$\theta(St(\tau_D) \boxtimes \mu) \ne 0$ iff $\chi$ is non-trivial
and $\tau_D \otimes \chi = \tau_D$, in which case
\[  \Theta(St(\tau_D) \boxtimes \mu)=  \theta(St(\tau_D) \boxtimes \mu) = St(\chi,  JL(\tau_D)), \]
which is generic. \vskip 5pt

(iii) ({\bf Generalized Speh}) Similarly, with $\mu = \omega_{\tau_D} \cdot \chi$, $\theta(Sp(\tau_D) \boxtimes \mu)$ is nonzero iff  $\tau_D \otimes \chi = \tau_D$. Suppose that
this holds. Then if $\chi \ne 1$,
\[  \Theta(Sp(\tau_D) \boxtimes \mu) =  \theta(Sp(\tau_D) \boxtimes \mu) =  Sp(\chi,  JL(\tau_D)), \]
which is non-generic. If $\chi =1$,
\[  \Theta(Sp(\tau_D) \boxtimes \omega_D)  = \theta(Sp(\tau_D) \boxtimes \omega_D)=
I_{Q(Z)}(|-|, JL(\tau_D) \cdot |-|^{-1/2}),   \] which is generic.
\vskip 10pt

(iv) ({\bf Twisted Steinberg})  If $St_{\chi}$ is the twisted Steinberg
representation of $GL_2(D)$, then $\theta(St_{\chi} \boxtimes \mu)
\ne 0$ iff $\mu = \chi^2$, in which case
\[  \theta(St_{\chi} \boxtimes \chi^2) = St_{PGSp_4} \otimes \chi. \]
\end{Thm}
\vskip 10pt

\begin{Thm}  \label{T:local-theta}
Consider the case when $D$ is split. \vskip 5pt

(i) ({\bf Principal series}) The irreducible principal series
representation $PS(\tau_1, \tau_2) \boxtimes \mu$, with $\tau_i$
discrete series representations  and $\mu^2 = \omega_1 \cdot
\omega_2$, participates in the local theta correspondence with
$GSp_4$ iff one of the following holds: \vskip 5pt

\begin{itemize}
\item[(a)]  $\omega_1 = \omega_2 = \mu$;

\item[(b)] $\tau_1^{\vee} \cong \tau_2 \otimes \mu^{-1}$ and  $\omega_1 \ne \mu$.
\end{itemize}

\vskip 5pt Moreover, if (a) holds, then $\theta(PS(\tau_1, \tau_2)
\boxtimes \mu)$ is the generic representation of $GSp_4(F)$ which is
the theta lift of the representation $\tau_1 \boxtimes \tau_2$ of
$GSO(D)$. If $\tau_1 = \tau_2 = \tau$, then $\theta(PS(\tau_1,
\tau_2) \boxtimes \mu)$ is the unique generic summand of the
tempered representation $I_{Q(Z)}(1, \tau)$, which is denoted by $\pi_{gen}({\tau})$ in \cite{GT}.
\vskip 5pt

If (b) holds,  then
\[  \theta(PS(\tau_1, \tau_2) \boxtimes \mu)=I_{Q(Z)}(\frac{\mu}{\omega_1}, \tau_2 \cdot \frac{\omega_1}{\mu})
= I_{Q(Z)}(\frac{\omega_1}{\mu}, \tau_2) \] which is
irreducible and generic.

\vskip 5pt
 (ii) ({\bf Generalized Steinberg})  Suppose that $\tau$ is supercuspidal with central character $\omega_{\tau}$ and $\mu = \omega_{\tau} \cdot \chi$ with $\chi^2 = 1$. Then $\theta(St(\tau) \boxtimes \mu) \ne 0$ iff $\chi$ is non-trivial and $\tau \otimes \chi = \tau$, in which
case
\[  \Theta(St(\tau) \boxtimes \mu) =  \theta(St(\tau) \boxtimes \mu) = St(\chi,  \tau), \]
which is generic. \vskip 5pt

(iii) ({\bf Generalized Speh}) Similarly, $\theta(Sp(\tau) \boxtimes \mu)$
is nonzero iff  $\tau \otimes \chi = \tau$. Suppose that this holds.
Then if $\chi \ne 1$,
\[   \Theta(Sp(\tau) \boxtimes \mu) =  \theta(Sp(\tau) \boxtimes \mu) =  Sp(\chi,  \tau), \]
which is non-generic. If $\chi =1$, then
\[  \Theta(Sp(\tau) \boxtimes \omega_{\tau})  = \theta(Sp(\tau) \boxtimes \omega_{\tau})=
I_{Q(Z)}(|-|, \tau \cdot |-|^{-1/2}), \]  which is generic.
\vskip 5pt

(iv) ({\bf Twisted Steinberg}) If $St_{\chi}$ is the twisted Steinberg
representation of $GL_4$, then $\theta(St_{\chi} \boxtimes \mu) \ne
0$ iff $\mu = \chi^2$, in which case
\[  \theta(St_{\chi} \boxtimes \chi^2) = St_{PGSp_4} \otimes \chi. \]

\end{Thm}
\vskip 10pt

Before coming to the proofs of the Theorems, let us draw a number of
consequences. Firstly, a comparison of Thm \ref{T:local-theta-D}(ii,
iv) and Thm. \ref{T:local-theta}(ii, iv) gives the purely local
proof of Prop. \ref{P:JL} for non-supercuspidal discrete series
representations promised at the end of the previous section. Indeed,
one has: \vskip 5pt

\begin{Cor}
(i) Let $\pi$ be a discrete series representation of $GL_4$ and
$\pi_D$ its Jacquet-Langlands lift to $GL_2(D)$. If $\pi = St(\tau)$
with $\tau$ supercuspidal,  then $\pi$ (resp. $\pi_D$) has Shalika
period with respect to $\mu$ iff $\mu = \omega_{\tau} \cdot \chi$
where $\chi$ is a non-trivial quadratic character such that $\tau
\otimes \chi \cong \tau$. When this holds, the (big $=$ small) theta lifts of $\pi$
and $\pi_D$ are isomorphic as representations of $GSp_4$.

\vskip 5pt

(ii) If $\pi = St_{\chi}$ is a twisted Steinberg representation,
then $\pi$ (resp. $\pi_D$) has Shalika period with respect to $\mu$
iff $\mu = \chi^2$, in which case the small theta lifts of $\pi$ and
$\pi_D$ are isomorphic as representations of $GSp_4$.
\end{Cor}
\vskip 10pt

Secondly, Thm. \ref{T:local-theta-D}(i) completes the proof of the
last statement in Thm. \ref{T:non-split}, regarding non-genericity
of the small theta lift.

\vskip 10pt Lastly, the theorems allow one to determine whether the
generalized Speh representations possess Shalika periods with
respect to $\mu$. This answers a question raised by Prasad in [P2].
\vskip 5pt

\begin{Thm}  \label{T:speh}
(i) The generalized Speh representation $Sp(\tau)$ has Shalika
period with respect to $\mu$ if and only if $\mu = \omega_{\tau}$.
\vskip 5pt

(ii) The generalized Speh representation $Sp(\tau_D)$ has Shalika
period with respect to $\mu$ if and only if $\mu = \omega_{\tau_D}$.
\end{Thm}

\begin{proof}
The proofs of (i) and (ii) are similar, so we shall only address
(i). By Thm. \ref{T:local-theta}(iii), we see that with $\mu =
\omega_{\tau} \cdot \chi$, $\Theta(Sp(\tau) \boxtimes \mu)$ is
generic if and only if $\chi =1$. Thus, $Sp(\tau)$ has Shalika
period with respect to $\mu$ if and only if  $\chi=1$.
\end{proof}

\vskip 15pt

We must now prove Thms. \ref{T:local-theta-D} and
\ref{T:local-theta}.  Given the essential similarity in the
statements of the two theorems, it is not surprising that one can
execute their proofs concurrently. Thus, in the remainder of the
section, $D$ is a possibly split quaternion algebra. \vskip 10pt

The key step is the computation of the normalized Jacquet module of
$\Omega_D$ with respect to $Q(Z)$ and $P(J)$. This is a
by-now-standard computation, following the lines of [K], and we
shall simply state the results below. For the computations, it is in
fact better not to identify $GSO(V_D)$ with a quotient of $GL_2(D)
\times GL_1$. Thus, we shall work directly with the parabolic $P(J)
= M(J) \cdot N(J)$ with $M(J) = GL(J) \times GSO(D)$, and we
represent an element of $M(J)$ by $(a, \alpha, \beta)$ with
\[ (\alpha,\beta) \in GSO(D) \cong (D^{\times} \times D^{\times})/\{(z,z^{-1}): z \in F^{\times} \}. \]
For a character $\chi$ and a representation $\tau_1 \boxtimes
\tau_2$ of $GSO(D)$, one may consider the normalized induced
representation $I_{P(J)}(\chi, \tau_1 \boxtimes \tau_2)$. \vskip 5pt

The relation of the two descriptions of principal series representations of $GSO(V_D)$ is as follows.
 Suppose that under the natural map $P_D \times GL_1 \longrightarrow P(J)$,
\[  \left( \left( \begin{array}{cc}
\alpha & \\
 & \beta  \end{array} \right), z \right)  \mapsto (a, \alpha' , \beta'), \]
 then we have:
 \[  \begin{cases}
 \alpha = a \cdot \overline{\alpha'}^{-1} \\
 \beta  = \beta' \\
 z = a^{-1}\cdot  \mathbb{N}(\alpha'). \end{cases} \]
 From this, one deduces that
 \[  PS(\tau_1,\tau_2)\boxtimes \mu  \cong I_{P(J)}( \omega_1\mu^{-1} , \tau_1^{\vee} \mu \boxtimes \tau_2). \]
   \vskip 10pt

Now we have: \vskip 10pt

\begin{Prop} \label{P:jacquet1}
Let $R_{P(J)}(\Omega_D)$ denote the normalized Jacquet module of
$\Omega_D$ along $P(J)$ (where $D$ is possibly split). Then we have
a short exact sequence of $M(J) \times GSp(W)$-modules: \vskip 5pt
\[  \begin{CD}
0 @>>> A @>>> R_{P(A)}(\Omega_D) @>>> B @>>> 0. \end{CD} \]
Here,
\[  B \cong \Omega_{W, D}, \]
\vskip 5pt \noindent where $\Omega_{W,D}$ is the induced Weil
representation for $GSp(W) \times GSO(D)$, and \vskip 3pt
\[  A \cong   I_{Q(Z)} \left(  S(F^{\times})   \otimes \Omega_{W',D}  \otimes   |\lambda_{W'}|^{-1}  \otimes  |\lambda_D|^{-1} \right). \]
\vskip 5pt \noindent  The action of $(GL(J) \times GSO(D))\times
(GL(Z) \times GSp(W'))$ on $S(F^{\times})$ is given by: \vskip 5pt
\[  ((a, h), (b, g)) \cdot f (x) = f(b^{-1}\cdot x \cdot a  \cdot \lambda_{W'}(g)), \]
\vskip 5pt \noindent and $\Omega_{W',D}$ denotes the induced Weil
representation of $GSp(W') \times GSO(D)$.
\end{Prop}
 \vskip 10pt

 \begin{Prop} \label{P:jacquet2}
 Let $R_{Q(Z)}(\Omega_D)$ denote the normalized Jacquet module of $\Omega_D$ along $Q(Z)$ (where $D$ is possibly split). Then we have a short exact sequence of
 $GSO(V_D) \times L(Z)$-modules:
 \vskip 5pt
 \[ \begin{CD}
 0 @>>> A'@>>> R_{Q(Z)}(\Omega_D) @>>> B' @>>> 0.\end{CD}  \]
 Here,
 \[  B' \cong |det_Z|  \boxtimes  \Omega_{W', V_D} \]
 \vskip 5pt
 \noindent where $\Omega_{W',V_D}$ is the induced Weil representation of $GSp(W') \times GSO(V_D)$ and
 \vskip 3pt
 \[  A' \cong  I_{P(J)} \left( S(F^{\times}) \otimes \Omega_{W', D}  \otimes |det_Z| \cdot |det_J|^{-1} \cdot
 |\lambda_{W'}|^{-2} \cdot |\lambda_D|^{-1}  \right). \]
\vskip 5pt \noindent The action of $(GL(J) \times GSO(D)) \times
(GL(Z) \times GSp(W'))$ on $S(F^{\times})$ is given by \vskip 5pt
\[  ((a,h), (b,g)) \cdot f(x) = f( a^{-1} \cdot \lambda_{W'}(g)^{-1} \cdot x \cdot b), \]
\vskip 5pt
\noindent  and $\Omega_{W',D}$ is the induced Weil representation of
$GSp(W') \times GSO(D)$.
 \end{Prop}
 \vskip 10pt

Applying Frobenius reciprocity and Props. \ref{P:jacquet1} and
\ref{P:jacquet2}, we obtain: \vskip 5pt

\begin{Prop}  \label{P:Hom}
(i)  Consider the space
\[  \Hom_{GSO(V_D)}(\Omega_D, I_{P(J)}(\chi, \tau_1 \boxtimes \tau_2) ) \]
as a representation of $GSp(W)$. Then we have: \vskip 5pt

\begin{itemize}
\item[(a)] If $\chi \ne 1$, then
\[   \Hom_{GSO(V_D)}(\Omega_D, I_{P(J)}(\chi, \tau_1 \boxtimes \tau_2))
 = 0 \]
unless
\[  \tau_1 = \tau_2 = \tau, \]
in which case
\[ \Hom_{GSO(V_D)}(\Omega_D, I_{P(J)}(\chi, \tau \boxtimes \tau)) =  I_{Q(Z)}\left( \chi^{-1}, JL(\tau) \cdot \chi \right)^* \quad \text{(full linear dual)}. \]
\vskip 5pt

\item[(b)] If $\chi=1$ but $\tau_1 \ne \tau_2$, then
\[  \Hom_{GSO(V_D)}(\Omega_D, I_{P(J)}(\chi, \tau_1 \boxtimes \tau_2) )
= \Theta_{W,D}(\tau_1 \boxtimes  \tau_2)^*, \] where
$\Theta_{W,D}(\tau_1 \boxtimes \tau_2)$ denotes the big theta lift
of $\tau_1 \boxtimes \tau_2$ from $GSO(D)$ to $GSp(W)$. \vskip 5pt

\item[(c)] If $\chi = 1$ and $\tau_1 = \tau_2 = \tau$, then we have an exact sequence:
\[  \begin{CD} 0 @>>> \Theta_{W,D}(\tau \boxtimes \tau)^*@>>>
\Hom_{GSO(V_D)}(\Omega_D, I_{P(J)}(\chi, \tau \boxtimes \tau)) @>>>
 \left( I_{Q(Z)}(1, JL(\tau)) \right)^*. \end{CD} \]
\end{itemize}
\vskip 10pt

(ii) Assume that $\chi \ne |-|$. Then as a representation of
$GSO(V_D)$,
\[  \Hom_{GSp(W)}(\Omega_D, I_{Q(Z)}(\chi, \tau)) = I_{P(J)}(\chi^{-1}, (JL(\tau) \cdot \chi) \boxtimes (JL(\tau) \cdot \chi))^*. \]
\end{Prop}
\vskip 10pt

\begin{proof}[Proof of Theorems \ref{T:local-theta-D} and \ref{T:local-theta}]
Now we can prove Thms. \ref{T:local-theta-D} and
\ref{T:local-theta}. In the following, we shall use the fact that if
$\pi$ is an irreducible representation of $GSO(V_D)$, then
\[  \Theta(\pi)^* \cong \Hom_{GSO(V_D)}(\Omega_D,  \pi). \]
We consider the different cases separately. \vskip 10pt

\noindent{\bf \underline{Principal Series}} \vskip 5pt

Suppose that
\[  PS(\tau_1,\tau_2)\boxtimes \mu  \cong I_{P(J)}( \omega_1\mu^{-1} , \tau_1^{\vee} \mu \boxtimes \tau_2) \]
is an irreducible principal series representation with $JL(\tau_i)$
discrete series representations. If $\omega_1 \ne \mu$, then by
Prop. \ref{P:Hom}(i)(a), we deduce that
\[  \Theta(PS(\tau_1,\tau_2)\boxtimes \mu) = 0 \]
unless $\tau_1^{\vee} \cdot \mu = \tau_2$, in which case
\[   \Theta(PS(\tau_1,\tau_2)\boxtimes \mu)= I_{Q(Z)}(\frac{\mu}{\omega_1}, JL(\tau_2) \cdot \frac{\omega_1}{\mu}). \]
Since the latter is irreducible also, it is isomorphic to
$I_{Q(Z)}(\omega_1 \mu^{-1}, JL(\tau_2))$. \vskip 5pt

On the other hand, suppose that $\omega_1 = \mu$ but $\tau_1 =
\tau_1^{\vee} \cdot \mu  \ne \tau_2$. Then Prop. \ref{P:Hom}(i)(b)
shows that
\[    \Theta(PS(\tau_1,\tau_2)\boxtimes \mu) = \Theta_{W, D}(\tau_1 \boxtimes \tau_2). \]
\vskip 5pt

Finally,  suppose that $\omega_1 = \mu$ and $\tau_1 = \tau_2 =
\tau$. By Prop. \ref{P:Hom}(i)(c), we obtain
\[ \begin{CD}
I_{Q(Z)}(1, JL(\tau)) @>>> \Theta(PS(\tau,\tau)\boxtimes \mu) @>>>
\Theta_{W,D}(\tau \boxtimes  \tau)@>>>0.
\end{CD} \]
This shows that
\[  \theta(PS(\tau,\tau)\boxtimes \mu) \supset \theta_{W,D}(\tau \boxtimes \tau). \]
We now have to examine if the two constituents of
$I_{Q(Z)}(1,JL(\tau))$ contribute to
 $\Theta(PS(\tau,\tau)\boxtimes \mu)$ or  even $\theta(PS(\tau,\tau)\boxtimes \mu)$.
\vskip 5pt

We shall suppose that $D$ is non-split; the split case is similar
and so we omit the details. Then $\theta_{W,D}(\tau \boxtimes \tau)$
is the unique non-generic summand of $I_{Q(Z)}(1, JL(\tau))$. Now by
Prop. \ref{P:prasad}(ii), one knows that $PS(\tau \boxtimes \tau)$
has Shalika period with respect to $\mu = \omega_1$.  Thus,
$\Theta(PS(\tau,\tau)\boxtimes \mu)$ is generic, so that the generic
summand of $I_{Q(Z)}(1,JL(\tau))$ does occur as a submodule of $
\Theta(PS(\tau,\tau)\boxtimes \mu)$. However, it does not occur as a
quotient of $ \Theta(PS(\tau,\tau)\boxtimes \mu)$. This follows
by [KR, Thm. 3.8]: since the generic summand of
$I_{Q(Z)}(1,JL(\tau))$ has nonzero theta lift to $GSO(2,2)$, it
cannot participate in the theta correspondence with $GSO(V_D)$.
Thus, we now know:
\[   \theta(PS(\tau,\tau)\boxtimes \mu) =   \theta_{W,D}(\tau \boxtimes \tau) \quad \text{or} \quad
2 \cdot  \theta_{W,D}(\tau \boxtimes \tau). \] \vskip 5pt

To show that the latter is not possible, we use Prop.
\ref{P:Hom}(ii) to see that
\[  \Hom_{GSp(W)}(\Omega_D, I_{Q(Z)}(1,JL(\tau)) = I_{P(J)}(1, \tau \boxtimes \tau)^*= (PS(\tau, \tau) \boxtimes  \omega_{\tau})^*  \]
so that
\[  PS(\tau,\tau) \boxtimes \omega_{\tau} \twoheadrightarrow  \Theta( \theta_{W,D}(\tau \boxtimes \tau)).
\]
This shows that \vskip 5pt
\[  \begin{cases}
 \theta(PS(\tau,\tau)\boxtimes \mu) =   \theta_{W,D}(\tau \boxtimes \tau);  \\
\Theta( \theta_{W,D}(\tau \boxtimes \tau)) = PS(\tau,\tau)\boxtimes
\mu.  \end{cases} \] \vskip 5pt \noindent  This completes the proof
of Thms. \ref{T:local-theta-D}(i) and \ref{T:local-theta}(i). \vskip
15pt

\noindent{\bf \underline{Generalized Steinberg and Speh}} \vskip 5pt

Now we consider the theta lift of the generalized Steinberg
representation $St(\tau) \otimes \mu$ and the generalized Speh
representation $Sp(\tau) \boxtimes \mu$, where
\[  \mu = \omega_{\tau} \cdot \chi \quad \text{with $\chi^2 =1$.} \]
Since
\[  St(\tau) \boxtimes \mu \hookrightarrow I_{P(J)}(\chi |-|, (\tau \cdot \chi |-|^{-1/2}) \boxtimes (\tau|-|^{-1/2})), \]
we deduce by Prop. \ref{P:Hom}(i)(a) that
\[  \Theta(St(\tau) \boxtimes \mu)^* \hookrightarrow \Hom_{GSO(V_D)}(\Omega_D, I_{P(J)}(\chi |-|, (\tau \cdot \chi |-|^{-1/2}) \boxtimes (\tau|-|^{-1/2})), \]
which vanishes unless $\tau \otimes \chi \cong \tau$,  in which case
one has:
\[  I_{Q(Z)}(\chi |-|^{-1},  JL(\tau) \cdot |-|^{1/2}) \twoheadrightarrow  \Theta(St(\tau) \boxtimes \mu).
\]
Recall that the above induced representation is irreducible if $\chi
=1$ and has $St(\chi , JL(\tau))$ as
unique irreducible quotient if $\chi \ne 1$. From this, we conclude
that if $\chi \ne 1$, then one has: \vskip 5pt
\begin{itemize}
\item  $\theta(St(\tau) \boxtimes \mu)\subset St(\chi, JL(\tau))$;
\vskip 5pt

\item $\theta(St(\chi, JL(\tau))) \ne 0$,
\end{itemize}
\vskip 5pt \noindent whereas if $\chi =1$, one has \vskip 5pt

\begin{itemize}
\item  $\theta(St(\tau) \boxtimes \mu)\subset I_{Q(Z)}(|-|, JL(\tau) \cdot |-|^{-1/2})$
\vskip 5pt

\item $\theta(I_{Q(Z)}(|-|, JL(\tau) \cdot |-|^{-1/2})) \ne 0$.
\end{itemize}
\vskip 10pt

On the other hand, if $\chi \ne 1$, one may apply  Prop.
\ref{P:Hom}(ii) to $I_{Q(Z)}(\chi |-|, JL(\tau) \cdot |-|^{-1/2})$
and arguing as above, one obtains: \vskip 5pt
\begin{itemize}
 \item $\theta(St(\chi, JL(\tau))) \subset St(\tau) \boxtimes \mu$;

\item  $\theta(St(\tau) \boxtimes \mu) \ne 0$.
\end{itemize}
\vskip 5pt

\noindent  Hence, we have shown that when $\chi \ne 1$, \vskip 5pt
\[   \begin{cases}
\theta(St(\tau) \boxtimes \mu) =  St(\chi, JL(\tau)); \\
\theta(St(\chi, JL(\tau))) =  St(\tau)
\boxtimes \mu. \end{cases} \] \vskip 10pt

Similarly, by applying Prop. \ref{P:Hom}(i)(a) to $I_{P(J)}(\chi
|-|^{-1}, (\tau \cdot \chi |-|^{-1/2}) \boxtimes (\tau|-|^{1/2}))$
and Prop. \ref{P:Hom}(ii) to $I_{Q(Z)}(\chi |-|^{-1}, JL(\tau) \cdot
|-|^{1/2})$, one deduces that
\[  \Theta(Sp(\tau) \boxtimes \mu) = 0 \]
unless $\tau \otimes \chi  = \tau$, in which case one has, if $\chi
\ne 1$,
\[  \begin{cases}
 \theta(Sp(\tau) \boxtimes \mu) =  Sp(\chi, JL(\tau)); \\
\theta(Sp(\chi, JL(\tau))) =  Sp(\tau)
\boxtimes \mu, \end{cases} \] whereas if $\chi = 1$,
\[  \begin{cases}
 \theta(Sp(\tau) \boxtimes \mu) =  I_{Q(Z)}(|-|, JL(\tau) \cdot |-|^{-1/2}); \\
\theta(I_{Q(Z)}(|-|, JL(\tau) \cdot |-|^{-1/2})) =  Sp(\tau)
\boxtimes \mu. \end{cases} \] This then implies that (for $\chi =1$)
\[  \Theta(St(\tau) \boxtimes \omega_{\tau}) = 0. \]
\vskip 5pt We have more or less completed the proof of Thm.
\ref{T:local-theta-D}(ii) and (iii), as well as Thm.
\ref{T:local-theta}(ii) and (iii), except that we still need to
check that the small theta lifts and the big theta lifts are equal.
\vskip 10pt

For that, we argue by contradiction. Suppose for example that (when
$\chi \ne 1$)
\[   \Theta(St(\tau) \boxtimes \mu)  \ne  \theta(St(\tau) \boxtimes \mu). \]
Then we must have
\[   \Theta(St(\tau) \boxtimes \mu) =  I_{Q(Z)}(\chi |-|^{-1},  JL(\tau) \cdot |-|^{1/2}). \]
This means that
\[  (St(\tau) \boxtimes \mu)^* \hookrightarrow
\Hom_{GSp(W)}(\Omega_D,  I_{Q(Z)}(\chi |-|^{-1},  JL(\tau) \cdot
|-|^{1/2})). \] But by Prop. \ref{P:Hom}(ii),
\[  \Hom_{GSp(W)}(\Omega_D,  I_{Q(Z)}(\chi |-|^{-1},  JL(\tau) \cdot |-|^{1/2}))
= I_{P(J)}(\chi |-|, \tau \cdot |-|^{-1/2})^*. \] Thus, we would
conclude that
\[  I_{P(J)}(\chi |-|, \tau \cdot |-|^{-1/2}) \twoheadrightarrow  St(\tau) \boxtimes \mu, \]
which is a contradiction. The other cases are treated similarly; we
omit the details. This completes the proof of the parts of Thms.
\ref{T:local-theta-D} and \ref{T:local-theta} pertaining to the
generalized Steinberg and generalized Speh representations. \vskip
15pt

\noindent{\bf \underline{Twisted Steinberg Representations}}
\vskip 5pt

Now consider the twisted Steinberg representation $St_{\chi}$ with
$\chi^4 = \mu^2$. We assume that $D$ is non-split, since the case
for split $D$ is similar. Since
\[  St_{\chi} \boxtimes \mu  \hookrightarrow I_{P(J)}(\frac{\chi^2}{\mu}, \frac{\mu}{\chi} |-|^{-1} \boxtimes \chi|-|^{-1}), \]
and
\[  St_{PGSp_4} \otimes \chi \hookrightarrow I_{Q(Z)}(|-|^2, JL(\chi) |-|^{-1}), \]
we may apply Prop. \ref{P:Hom}(i) and (ii)  to conclude that a
necessary condition for the non-vanishing of theta lifts is $\mu =
\chi^2$, in which case a similar argument as the above cases shows
that
\[  \theta(St_{\chi} \boxtimes \chi^2) = St_{PGSp_4} \otimes \chi \]
and
\[   \theta( St_{PGSp_4} \otimes \chi ) = St_{\chi} \boxtimes \chi^2. \]
This completes the proof of Thms. \ref{T:local-theta-D}(iv) and
\ref{T:local-theta}(iv).
\end{proof}
\vskip 15pt

\noindent{\bf \underline{Explicit Theta Correspondence for $GL_2 \times GSO(V_D)$}}
\vskip 5pt

We conclude this section by describing the local theta correspondence for $GL_2 \times GSO(V_D)$, where $D$ is possibly split. This was needed at certain places in Section \ref{S:cde}.  Hence, let
$\Omega_{W',D}$ be the induced Weil representation for this dual pair which can be realized on
$S(V_D)$. As in Section \ref{S:cde}, we have $W' = F \cdot e \oplus F \cdot f$ and we let $B = T \cdot U$ be the Borel subgroup stabilizing $F \cdot e$. Then we have
\vskip 10pt

\begin{Prop}
Let $R_B(\Omega_{W',D})$ be the normalized Jacquet module of $\Omega_{W',D}$
with respect to the unipotent radical $U$ of $B$. There is a short exact sequence of representations of $T  \times GSO(V_D)$:
\[  \begin{CD}
0 @>>> I_{P(J)}( S(F^{\times} \times F^{\times})) @>>> R_B(\Omega_{W',D}) @>>> S(F^{\times}) @>>> 0.
\end{CD} \]
Here, the actions of $t(a,b) \in T$ and $h \in GSO(V_D)$  on $S(F^{\times})$ are given by
\[  \begin{cases}
(t(a,b)\cdot \phi)(t) = |a|^{-1/2} \cdot |b|^{-5/2} \cdot  \phi(tab) \\
(h \cdot \phi)(t) = |\lambda(h)|^{-3/2} \cdot \phi(t \lambda(h)). \end{cases} \]
On the other hand, the action of $T \times M(J) = T \times GL(J) \times GSO(D)$ on
$S(F^{\times} \times F^{\times})$ is given as follows. For $t(a,b) \in T$,
\[  (t(a,b)\cdot \phi)(t,x) = |a|^{-1/2} \cdot |b|^{-5/2} \cdot \phi(tab, b^{-1}x). \]
For $\alpha \in GL(J)$,
\[  (\alpha \cdot \phi)(t,x) = |\alpha|^{-2} \cdot \phi(t, \alpha^{-1}x), \]
and for $h \in GSO(D)$,
\[  (h \cdot \phi)(t,x) = |\lambda_D(h)|^{-1/2} \cdot \phi(t \lambda_D(h),x). \]
\end{Prop}
\vskip 10pt

Using this proposition, we deduce (at least for parts (i) and (ii)):
\vskip 10pt

\begin{Thm} \label{T:explicit}
Let $\tau$ be an irreducible infinite dimensional unitary (up to twisting) representation of $GL_2$.
\vskip 5pt

(i) If $\tau = \pi(\chi_1, \chi_2)$, then
\[  \Theta_D(\tau) = PS(\chi_1\circ \det , \chi_2\circ \det) \boxtimes (\chi_1 \chi_2) \]
which is irreducible.
\vskip 5pt

(ii) If $\tau = st_{\chi}$ is a twisted Steinberg representation, then when $D$ is non-split,
\[  \Theta_D(\tau) = PS(JL(\tau)|-|^{1/2}, JL(\tau)|-|^{-1/2}) \boxtimes \omega_{\tau} \]
which is irreducible. On the other hand, when $D$ is split, then $\theta_D(\tau)$ is the unique irreducible quotient
of $PS(\tau|-|^{1/2}, \tau|-|^{-1/2})$.
\vskip 5pt

(iii) If $\tau$ is supercuspidal, then $\theta_D(\tau)  = Sp(JL(\tau)) \boxtimes \omega_{\tau}$.
\end{Thm}
\vskip 10pt

In the interest of space and time, we leave the details of the proof to the reader.

\newpage

\begin{center}

\begin{table}
\centering \caption{Explicit theta lifts from $GSO(V_D)$ to $GSp_4$}
\label{maintable} \vspace{-3ex}
$$
\renewcommand{\arraystretch}{1.5}
 \begin{array}{|c|c|c|c|c|c|}
  \hline
   \multicolumn{4}{|c|}{\mbox{$\pi\boxtimes\mu\in\Irr(GSO(V_D))$}}
   &\multicolumn{2}{|c|}{\mbox{$\theta(\pi\boxtimes\mu)\in\Irr(GSp_4)$}}\\\hline\hline

   \multirow{4}{*}{$PS(\tau_{D,1},\tau_{D,2})\boxtimes\mu$}&\mbox{a}&\multirow{2}{*}{$\omega_{D,1}=\omega_{D,2}=\mu$}&\tau_{D,1}\neq\tau_{D,2}
   &\multirow{2}{*}{$\theta(\tau_{D,1}\boxtimes\tau_{D,2})$}
   &\mbox{non-generic S.C.}
   \\ \cline{4-4}\cline{6-6}
   &\mbox{b}&&\tau_{D,1}=\tau_{D,2}=\tau_D
   &&\pi_{ng}(JL(\tau_D))
   \\  \cline{3-6}
   &\mbox{c}&\multicolumn{2}{|c|}{\tau_{D,1}^\vee\cong\tau_{D,2}\otimes\mu^{-1}, \omega_{D,1}\neq\mu}
   &\multicolumn{2}{|c|}{I_{Q(Z)}(\omega_{D,1}\cdot\mu^{-1}, JL(\tau_{D,2}))}
   \\ \cline{3-6}
   &\mbox{d}&\multicolumn{2}{|c|}{\mbox{otherwise}}&\multicolumn{2}{|c|}{$0$}
   \\
   \hline

   St(\tau_D)\boxtimes\mu, \dim\tau_D>1&\mbox{a}&\multicolumn{2}{|c|}{\chi\neq1, \tau_D\otimes\chi=\tau_D}
   &\multicolumn{2}{|c|}{St(\chi,JL(\tau_D))}
   \\ \cline{3-6}
   \mu=\omega_D\cdot\chi, \,  \, \chi^2=1&\mbox{b}&\multicolumn{2}{|c|}{\mbox{otherwise}}
   &\multicolumn{2}{|c|}{$0$}
   \\
   \hline

   \multirow{2}{*}{$Sp(\tau_D)\boxtimes\mu, \dim\tau_D>1$}
   &\mbox{a}&\multirow{2}{*}{$\tau_D\otimes\chi=\tau_D$}&\chi\neq 1
   &\multicolumn{2}{|c|}{Sp(\chi,JL(\tau_D))}
   \\ \cline{4-6}
   \multirow{2}{*}{$\mu=\omega_D\cdot\chi, \, \, \chi^2=1$}&\mbox{b}&&\chi=1
   &\multicolumn{2}{|c|}{I_{Q(Z)}(|-|,JL(\tau_D)\cdot |-|^{-1/2})}
   \\ \cline{3-6}
   &\mbox{c}&\multicolumn{2}{|c|}{\mbox{otherwise}}&\multicolumn{2}{|c|}{0}
   \\
   \hline

   \multirow{2}{*}{$St_{\chi}\boxtimes\mu$}&\mbox{a}&\multicolumn{2}{|c|}{\mu=\chi^2}
   &\multicolumn{2}{|c|}{St_{PGSp_4\otimes\chi}}
   \\ \cline{3-6}
   &\mbox{b}&\multicolumn{2}{|c|}{\mbox{otherwise}}
   &\multicolumn{2}{|c|}{0}
   \\
   \hline
\end{array}
$$

\end{table}

\end{center}

\begin{center}


\begin{table}[h]
\centering \caption{Explicit theta lifts from $GSO(V)$ to $GSp_4$}
\vspace{-3ex}
$$
\renewcommand{\arraystretch}{1.5}
 \begin{array}{|c|c|c|c|c|c|}
  \hline
   \multicolumn{4}{|c|}{\mbox{$\pi\boxtimes\mu\in\Irr(GSO(V))$}}
   &\multicolumn{2}{|c|}{\mbox{$\theta(\pi\boxtimes\mu)\in\Irr(GSp_4)$}}\\\hline\hline

   &\mbox{a}&\multirow{2}{*}{$\omega_{1}=\omega_{2}=\mu$}&\tau_{1}\neq\tau_{2}
   &\multirow{2}{*}{$\theta(\tau_{1}\boxtimes\tau_{2})$}
   &\mbox{generic}
   \\ \cline{4-4}\cline{6-6}
   PS(\tau_{1},\tau_{2})\boxtimes\mu&\mbox{b}&&\tau_{1}=\tau_{2}=\tau
   &&\pi_{gen}(\tau)
   \\  \cline{3-6}

   \mbox{$\tau_1$, $\tau_2$ discrete series}&\mbox{c}
   &\multicolumn{2}{|c|}{\tau_{1}^\vee\cong\tau_{2}\otimes\mu^{-1}, \omega_{1}\neq\mu}
   &\multicolumn{2}{|c|}{I_{Q(Z)}(\omega_{1}\cdot\mu^{-1}, \tau_{2})}
   \\ \cline{3-6}
   &\mbox{d}&\multicolumn{2}{|c|}{\mbox{otherwise}}&\multicolumn{2}{|c|}{$0$}
   \\
   \hline

   \mbox{$St(\tau)\boxtimes\mu,\ \tau$ supercuspidal}&\mbox{a}&\multicolumn{2}{|c|}{\chi\neq1, \tau\otimes\chi=\tau}
   &\multicolumn{2}{|c|}{St(\chi,\tau)}
   \\ \cline{3-6}
   \mu=\omega\cdot\chi, \, \, \chi^2=1&\mbox{b}&\multicolumn{2}{|c|}{\mbox{otherwise}}
   &\multicolumn{2}{|c|}{$0$}
   \\
   \hline

   \multirow{2}{*}{$Sp(\tau)\boxtimes\mu$, $\tau$ supercuspidal}&\mbox{a}&\multirow{2}{*}{$\tau\otimes\chi=\tau$}&\chi\neq 1
   &\multicolumn{2}{|c|}{Sp(\chi,\tau)}
   \\ \cline{4-6}
   \multirow{2}{*}{$\mu=\omega_D\cdot\chi, \, \, \chi^2=1$}&\mbox{b}&&\chi=1
   &\multicolumn{2}{|c|}{I_{Q(Z)}(|-|,\tau\cdot |-|^{-1/2})}
   \\ \cline{3-6}
   &\mbox{c}&\multicolumn{2}{|c|}{\mbox{otherwise}}&\multicolumn{2}{|c|}{0}
   \\
   \hline

   \multirow{2}{*}{$St_{\chi}\boxtimes\mu$}&\mbox{a}&\multicolumn{2}{|c|}{\mu=\chi^2}
   &\multicolumn{2}{|c|}{St_{PGSp_4\otimes\chi}}
   \\ \cline{3-6}
   &\mbox{b}&\multicolumn{2}{|c|}{\mbox{otherwise}}
   &\multicolumn{2}{|c|}{0}
   \\
   \hline
\end{array}
$$

\end{table}
\end{center}

\begin{center}
\begin{table}[h]
\centering \caption{Explicit theta lifts from $GL_2$ to $GSO(V_D)$}
\vspace{-3ex}
$$
\renewcommand{\arraystretch}{1.5}
 \begin{array}{|c|c|c|}
  \hline
   \multicolumn{2}{|c|}{\tau\in\Irr(GL_2)}
   &\theta(\tau)\in\Irr(GSO(V_D))\\\hline\hline
    \multicolumn{2}{|c|}{\pi(\chi_1,\chi_2)}&PS(\chi_1\circ\det, \chi_2\circ\det)\boxtimes(\chi_1\chi_2)
    \\ \hline
    \multirow{2}{*}{$st_\chi$}&\mbox{$D$ non-split}&PS(JL(\tau)|-|^{1/2}, JL(\tau)|-|^{-1/2})\boxtimes\omega_\tau
    \\ \cline{2-3}
    &\mbox{$D$ split}&\mbox{unique quotient of $PS(\tau|-|^{1/2}, \tau|-|^{-1/2})\boxtimes\omega_\tau$}
    \\ \hline
    \multicolumn{2}{|c|}{\mbox{Supercuspidal}}&Sp(JL(\tau))\boxtimes\omega_\tau
    \\ \hline

\end{array}
$$
\end{table}
\end{center}

\eject

\end{document}